# A recursive online algorithm for the estimation of time-varying ARCH parameters

RAINER DAHLHAUS[*] and SUHASINI SUBBA RAO[**]

*Institut für Angewandte Mathematik, Universität Heidelberg, Im Neuenheimer Feld 294, 69120 Heidelberg, Germany. E-mail: [*]dahlhaus@statlab.uni-heidelberg.de; [**]suhasini@stat.tamu.edu*

In this paper we propose a recursive online algorithm for estimating the parameters of a time-varying ARCH process. The estimation is done by updating the estimator at time point $t-1$ with observations about the time point $t$ to yield an estimator of the parameter at time point $t$. The sampling properties of this estimator are studied in a non-stationary context – in particular, asymptotic normality and an expression for the bias due to non-stationarity are established. By running two recursive online algorithms in parallel with different step-sizes and taking a linear combination of the estimators, the rate of convergence can be improved for parameter curves from Hölder classes of order between 1 and 2.

## 1. Introduction

The class of autoregressive conditional heteroscedastic (ARCH) processes can be generalized to include non-stationary processes, by including models with parameters which are time-dependent. More precisely, $\{X_{t,N}\}$ is called a time-varying ARCH (tvARCH) process of order $p$ if it satisfies the representation

$$X_{t,N} = Z_t \sigma_{t,N}, \qquad \sigma_{t,N}^2 = a_0\left(\frac{t}{N}\right) + \sum_{j=1}^{p} a_j\left(\frac{t}{N}\right) X_{t-j,N}^2, \tag{1}$$

where $\{Z_t\}$ are independent, identically distributed random variables with $\mathbb{E}(Z_0) = 0$ and $\mathbb{E}(Z_0^2) = 1$. This class of tvARCH processes was investigated in Dahlhaus and Subba Rao [4]. It was shown that it can locally be approximated by a stationary process; we summarize the details below. Furthermore, a local quasi-likelihood method was proposed to estimate the parameters of the tvARCH($p$) model.

A potential application of the tvARCH process is to model long financial time series. The modelling of financial data using non-stationary time series models has recently attracted considerable attention. A justification for using such models can be found, for example, in Mikosch and Stărică [9, 10]. However, given that financial time series are often sampled at high frequency, evaluating the likelihood as each observation comes online can be computationally expensive. Thus an 'online' method, which uses the previous estimate of the parameter at time point $t-1$ and the observation at time point $t$ to estimate the parameter at time point $t$ would be ideal and cost-effective. There exists a huge literature on recursive algorithms, mainly in the context of linear systems (cf. Ljung and Söderström [8]; Solo [12, 13]) or neural networks (cf. White [15]; Chen and White [2]). For a general overview, see also Kushner and Yin [7]. Motivated by the least mean squares algorithm in Moulines *et al.* [11] for time-varying autoregressive processes, we consider in this paper the following recursive online algorithm for tvARCH models:

$$\underline{\hat{a}}_{t,N} = \underline{\hat{a}}_{t-1,N} + \lambda \{X_{t,N}^2 - \underline{\hat{a}}_{t-1,N}^\top \mathcal{X}_{t-1,N}\} \frac{\mathcal{X}_{t-1,N}}{|\mathcal{X}_{t-1,N}|_1^2}, \qquad t = (p+1), \dots, N, \qquad (2)$$

with $\mathcal{X}_{t-1,N}^\top = (1, X_{t-1,N}^2, \dots, X_{t-p,N}^2)$, $|\mathcal{X}_{t-1,N}|_1 = 1 + \sum_{j=1}^p X_{t-j,N}^2$ and initial conditions $\underline{\hat{a}}_{p,N} = (0, \dots, 0)$. This algorithm is linear in the estimators, despite the nonlinearity of the tvARCH process. We call the stochastic algorithm defined in (2) the ARCH normalized recursive estimation (ANRE) algorithm. Let $\underline{a}(u)^\top = (a_0(u), \dots, a_p(u))$; then $\underline{\hat{a}}_{t,N}$ is regarded as an estimator of $\underline{a}(t/N)$ or of $\underline{a}(u)$ if $|t/N - u| < 1/N$.

In this paper we will prove the consistency and asymptotic normality of this recursive estimator. Furthermore, we will discuss the improvements of the estimator obtained by combining two estimates from (2) with different $\lambda$. Unlike in most other work in the area of recursive estimation the properties of the estimator are proved under the assumption that the true process is a process with time-varying coefficients, that is, a non-stationary process. The rescaling of the coefficients in (1) to the unit interval corresponds to the 'infill asymptotics' in nonparametric regression: as $N \to \infty$ the system does not describe the asymptotic behaviour of the system in a physical sense, but is meant as a meaningful asymptotics to approximate, for example, the distribution of estimates based on a finite sample size. A similar approach was used in Moulines *et al.* [11] for time-varying autoregressive models. A more detailed discussion of the relevance of this approach and the relation to non-rescaled processes can be found in Section 3.

In fact the ANRE algorithm resembles the NLMS algorithm investigated in Moulines *et al.* [11]. Rewriting (2), we have

$$\underline{\hat{a}}_{t,N} = \left(I - \lambda \frac{\mathcal{X}_{t-1,N} \mathcal{X}_{t-1,N}^\top}{|\mathcal{X}_{t-1,N}|_1^2}\right) \underline{\hat{a}}_{t-1,N} + \lambda \frac{X_{t,N}^2 \mathcal{X}_{t-1,N}}{|\mathcal{X}_{t-1,N}|_1^2}. \qquad (3)$$

We can see from (3) that the convergence of the ANRE algorithm relies on showing some type of exponential decay of the past. In this paper we will show that for any $p > 0$,

$$\mathbb{E}\left\|\prod_{i=1}^k \left(I - \lambda \frac{1}{|\mathcal{X}_{t-i-1,N}|_1^2} \mathcal{X}_{t-i-1,N} \mathcal{X}_{t-i-1,N}^\top\right)\right\|^p \le K(1 - \lambda\delta)^k \qquad \text{for some } \delta > 0, \quad (4)$$

where $\|\cdot\|$ denotes the spectral norm and $\prod_{i=0}^n A_i = A_0 \cdots A_n$. Roughly speaking, this means we have, on average, exponential decay of the past. Similar properties are often established in the control literature and referred to as persistence of excitation of the stochastic matrix (see, for example, Guo [5]; Aguech *et al.* [1]), which in our case is the matrix $(1 - \frac{\lambda}{|\mathcal{X}_{t-1,N}|_2^2} \mathcal{X}_{t-1,N} \mathcal{X}_{t-1,N}^\top)$. Persistence of excitation guarantees convergence of the algorithm, which we will use to prove the asymptotic properties of $\hat{\underline{a}}_{t,N}$.

In Section 2 we state all results on the asymptotic behaviour of $\hat{\underline{a}}_{t,N}$ including consistency, asymptotic normality and rate efficiency. Furthermore, we suggest a modified algorithm based on two parallel algorithms. In Section 3 we discuss practical implications. Sections 4 and 5 contain the proofs, which in large part are based on the perturbation technique. Some technical methods have been gathered in the Appendix. We note that some of results in the Appendix are of independent interest, as they deal with the probabilistic properties of ARCH and tvARCH processes and their vector representations.

# 2. The ANRE algorithm

We first review some properties of the tvARCH process. Dahlhaus and Subba Rao [4] and Subba Rao [14] have shown that the tvARCH process can be locally approximated by a stationary ARCH process. Let $u$ be fixed and

$$X_t(u) = Z_t \tilde{\sigma}_t(u), \qquad \sigma_t(u)^2 = a_0(u) + \sum_{j=1}^p a_j(u) X_{t-j}(u)^2, \tag{5}$$

where $\{Z_t\}$ are independent, identically distributed random variables with $\mathbb{E}(Z_0) = 0$ and $\mathbb{E}(Z_0^2) = 1$. We also set $\mathcal{X}_t(u)^\top = (1, X_t(u)^2, \ldots, X_{t-p+1}(u)^2)$. In Lemma 4.1 we show that $X_t(u)^2$ can be regarded as the stationary approximation of $X_{t,N}^2$ around the time points $t/N \approx u$.

**Assumption 2.1.** *Let $\{X_{t,N}\}$ and $\{X_t(u)\}$ be sequences of stochastic processes which satisfy* (1) *and* (5) *respectively.*

(i) *For some $r \in [1, \infty)$, there exists $\eta > 0$ such that $\{\mathbb{E}(Z_0^{2r})\}^{1/r} \sup_u \{\sum_{j=1}^p a_j(u)\} < 1 - \eta$.*

(ii) *There exists $0 < \rho_1 \le \rho_2 < \infty$ such that, for all $u \in (0, 1]$, $\rho_1 \le a_0(u) \le \rho_2$.*

(iii) *There exists $\beta \in (0, 1]$ and a constant $K$ such that for $u, v \in (0, 1]$,*

$$|a_j(u) - a_j(v)| \le K|u - v|^\beta \qquad for\ j = 0, \ldots, p.$$

(iv) *Let $Y_0(u) = a_0(u) Z_0^2$ and $Y_t(u) = \{a_0(u) + \sum_{j=1}^t a_j(u) Y_{t-j}(u)\} Z_t^2$ $(t = 1, \ldots, p)$. Define $\underline{Y}_p(u) = (1, Y_1(u), \ldots, Y_p(u))^\top$ and $\Sigma(u) = \mathbb{E}(\underline{Y}_p(u) \underline{Y}_p(u)^\top)$. Then there exists a constant $C$ such that $\inf_u \lambda_{\min}\{\Sigma(u)\} \ge C$.*

**Remark 2.1.** It is clear that $\Sigma(u)$ is a positive semi-definite matrix, hence its smallest eigenvalue is greater than or equal to zero. It can be shown that if $p/\mathbb{E}(Z_t^4)^{1/2} < 1$ and $\sup_u a_0(u) > 0$, then $\lambda_{\min}(\Sigma(u)) > (1 - p/\mathbb{E}(Z_t^4)^{1/2})/a_0(u)^{(2p+1)/(p+1)}$. However, this condition is only sufficient and lower bounds can be obtained under much weaker conditions.

We now investigate the asymptotic properties of the ANRE algorithm. We assume that $\lambda \to 0$ and $\lambda N \to \infty$ as $N \to \infty$. We mention explicitly that $\lambda$ does not depend on $t$, that is, we are considering the fixed-step-size case. The assumption $\lambda \to 0$ is possible in the triangle array framework of our model and the resulting assertions (e.g., Theorem 2.2) are meant as an approximation of the corresponding finite sample size distributions and not as the limit in any physical sense.

The following results are based on a representation proved at the end of Section 4.3. The difference $\hat{\underline{a}}_{t_0,N} - \underline{a}(u_0)$ is dominated by two terms, that is,

$$\hat{\underline{a}}_{t_0,N} - \underline{a}(u_0) = \mathcal{L}_{t_0}(u_0) + \mathcal{R}_{t_0,N}(u_0) + O_p(\delta_N), \tag{6}$$

where

$$\delta_N = \left( \frac{1}{(N\lambda)^{2\beta}} + \frac{\sqrt{\lambda}}{(N\lambda)^{\beta}} + \lambda + \frac{1}{N^{\beta}} \right), \tag{7}$$

$$\mathcal{L}_{t_0}(u_0) = \sum_{k=0}^{t_0-p-1} \lambda \{I - \lambda F(u_0)\}^k \mathcal{M}_{t_0-k}(u_0), \tag{8}$$

$$\mathcal{R}_{t_0,N}(u_0) = \sum_{k=0}^{t_0-p-1} \lambda \{I - \lambda F(u_0)\}^k \left( \left\{ \mathcal{M}_{t_0-k}\left( \frac{t_0-k}{N} \right) - \mathcal{M}_{t_0-k}(u_0) \right\} \right.$$
$$\left. + F(u_0) \left\{ \underline{a}\left( \frac{t_0-k}{N} \right) - \underline{a}(u_0) \right\} \right),$$

with

$$\mathcal{M}_t(u) = (Z_t^2 - 1)\sigma_t(u)^2 \frac{\mathcal{X}_{t-1}(u)}{|\mathcal{X}_{t-1}(u)|_1^2} \quad \text{and} \quad F(u) = \mathbb{E}\left( \frac{\mathcal{X}_0(u)\mathcal{X}_0(u)^\top}{|\mathcal{X}_0(u)|_1^2} \right). \tag{9}$$

We note that $\mathcal{L}_{t_0}(u_0)$ and $\mathcal{R}_{t_0,N}(u_0)$ play two different roles. $\mathcal{L}_{t_0}(u_0)$ is the weighted sum of the stationary random variables $\{X_t(u_0)\}_t$, which locally approximate the tvARCH process $\{X_{t,N}\}_t$, whereas $\mathcal{R}_{t_0,N}(u_0)$ is the (stochastic) bias due to non-stationarity; if the tvARCH process were stationary this term would be zero. It is clear from the above that the magnitude of $\mathcal{R}_{t_0,N}(u_0)$ depends on the regularity of the time-varying parameters $\underline{a}(u)$, for example, the Hölder class that $\underline{a}(u)$ belongs to. By using (6) we are able to obtain a bound for the mean squared error of $\hat{\underline{a}}_{t_0,N}$. Let $|\cdot|$ denote the Euclidean norm of a vector.

**Theorem 2.1.** *Suppose Assumption 2.1 holds with $r > 4$ and $u_0 > 0$. Then if $|u_0 - t_0/N| < 1/N$, we have*

$$\mathbb{E}\{|\hat{\underline{a}}_{t_0,N} - \underline{a}(u_0)|^2\} = O\left( \lambda + \frac{1}{(N\lambda)^{2\beta}} \right), \tag{10}$$

*where $\lambda \to 0$ as $N \to \infty$ and $N\lambda \gg (\log N)^{1+\varepsilon}$, with some $\varepsilon > 0$.*

The proof can be found at the end of Section 4.2. Theorem 2.1 implies $\hat{\underline{a}}_{t_0,N} \xrightarrow{\mathcal{P}} \underline{a}(u_0)$.

The stochastic term $\mathcal{L}_{t_0}(u_0)$ is the sum of martingale differences, which allows us to obtain the following central limit theorem, whose proof is at the end of Section 4.3.

**Theorem 2.2.** *Suppose Assumption* 2.1 *holds with* $r > 4$ *and* $u_0 > 0$. *Let* $\mathcal{R}_{t_0,N}(u_0)$ *be defined as in* (8). *If* $|t_0/N - u_0| < 1/N$

(i) *and* $\lambda \gg N^{-4\beta/(4\beta+1)}$ *and* $\lambda \gg N^{-2\beta}$, *then*

$$\lambda^{-1/2}\{\hat{\underline{a}}_{t_0,N} - \underline{a}(u_0)\} - \lambda^{-1/2}\mathcal{R}_{t_0,N}(u_0) \xrightarrow{\mathcal{D}} \mathcal{N}(0, \Sigma(u_0)); \tag{11}$$

(ii) *and* $\lambda \gg N^{-2\beta/(2\beta+1)}$ *then*

$$\lambda^{-1/2}\{\hat{\underline{a}}_{t_0,N} - \underline{a}(u_0)\} \xrightarrow{\mathcal{D}} \mathcal{N}(0, \Sigma(u_0)), \tag{12}$$

*where* $\lambda \to 0$ *as* $N \to \infty$ *and* $N\lambda \gg (\log N)^{1+\varepsilon}$, *for some* $\varepsilon > 0$, *with*

$$\Sigma(u) = \frac{\mu_4}{2} F(u)^{-1} \mathbb{E}\left(\frac{\sigma_1(u)^4 \mathcal{X}_0(u)\mathcal{X}_0(u)^\top}{|\mathcal{X}_0(u)|_1^4}\right), \qquad \mu_4 = \mathbb{E}(Z_0^4) - 1. \tag{13}$$

Until now we have assumed $\underline{a}(u) \in \mathrm{Lip}(\beta)$, where $\beta \leq 1$. Let $\dot{f}(u)$ denote the derivative of the vector or matrix $f(\cdot)$ with respect to $u$. Suppose $0 < \beta' \leq 1$ and $\dot{\underline{a}}(u) \in \mathrm{Lip}(\beta')$; then we say $\underline{a}(u) \in \mathrm{Lip}(1+\beta')$. We now show that an exact expression for the bias can be obtained if $\underline{a}(u) \in \mathrm{Lip}(1+\beta')$ and $\beta' > 0$. We make the following assumptions.

**Assumption 2.2.** *Let* $\{X_{t,N}\}$ *be a sequence of stochastic processes which satisfies Assumption* 2.1 *and* $|\dot{a}_i(u) - \dot{a}_i(v)| \leq K|u-v|^{\beta'}$ $(i = 0, \ldots, p)$ *for some* $\beta' > 0$.

Under this assumption we show in Lemma 5.3 that

$$\mathbb{E}\{\hat{\underline{a}}_{t_0,N} - \underline{a}(u_0)\} = -\frac{1}{N\lambda} F(u_0)^{-1} \dot{\underline{a}}(u_0) + O\left(\frac{1}{(N\lambda)^{1+\beta'}}\right). \tag{14}$$

We note that typically it is not possible to obtain an exact expression for the bias of parameter estimates of an ARCH process. By using the expression above for the bias we obtain the following theorem, whose proof is at the end of Section 5. Let $\mathrm{tr}(\cdot)$ denote the trace of the matrix.

**Theorem 2.3.** *Suppose Assumption* 2.2 *holds with* $r > 4$ *and* $u_0 > 0$. *Then if* $|t_0/N - u_0| < 1/N$, *we have*

$$\mathbb{E}|\hat{\underline{a}}_{t_0,N} - \underline{a}(u_0)|^2 = \lambda \, \mathrm{tr}\{\Sigma(u_0)\} + \frac{1}{(N\lambda)^2}|F(u_0)^{-1}\dot{\underline{a}}(u_0)|^2$$

$$+ O\left\{\frac{1}{(N\lambda)^{2+\beta'}} + \frac{\lambda^{1/2}}{(N\lambda)^{1+\beta'}} + \frac{1}{(N\lambda)^2}\right\}, \tag{15}$$

and if $\lambda$ is such that $\lambda^{-1/2}/(N\lambda)^{1+\beta'} \to 0$, then

$$\lambda^{-1/2}(\hat{\underline{a}}_{t_0,N} - \underline{a}(u_0)) + \lambda^{-1/2}\frac{1}{N\lambda}F(u_0)^{-1}\dot{\underline{a}}(u_0) \xrightarrow{\mathcal{D}} \mathcal{N}(0, \Sigma(u_0)), \qquad (16)$$

where $\lambda \to 0$ as $N \to \infty$ and $\lambda N \gg (\log N)^{1+\varepsilon}$, for some $\varepsilon > 0$.

Let $\|\underline{f}\|_\beta$ be the bounded Lipschitz norm and

$$\mathcal{C}^\beta(L) = \left\{ \underline{f}(\cdot) = (f_0(\cdot), \ldots, f_{p+1}(\cdot)) \Big| \underline{f} : [0,1] \to (\mathbb{R}^+)^{p+1}, \|\underline{f}\|_\beta \le L, \right.$$

$$\left. \sup_u \sum_{i=1}^{p+1} f_i(u) < 1 - \delta, 0 < \rho_1 < \inf_u f_0(u) < \rho_2 < \infty \right\}. \qquad (17)$$

In Dahlhaus and Subba Rao [3] we derived the following minimax risk for estimators $\hat{\underline{a}}_{t_0,N}$ of $\underline{a}(u_0)$:

$$\min_{\hat{\underline{a}}_{t_0,N} \in \sigma(X_{1,N}, \ldots, X_{N,N})} \max_{\underline{a}(u_0) \in \mathcal{C}^\nu(L)} \mathbb{E}|\hat{\underline{a}}_{t_0,N} - \underline{a}(u_0)|^2 \ge K N^{-2\nu/(2\nu+1)}. \qquad (18)$$

Comparing this bound with (10), it is straightforward to show that the ANRE algorithm attains the optimal rate if $\underline{a}(\cdot) \in \text{Lip}(\nu)$ with $\nu \le 1$ (with $\lambda = N^{-2\nu/(1+2\nu)}$). It is a different story when $1 < \nu < 2$. If $\underline{a}(\cdot) \in \text{Lip}(1 + \beta')$, $\beta' > 0$, the mean squared error of the ANRE estimator in (15) becomes minimal for $\lambda \approx N^{-2/3}$ with minimum rate $\mathbb{E}|\hat{\underline{a}}_{t_0,N} - \underline{a}(u_0)|^2 = O(N^{-2/3})$. However, in (18) the minimax rate in this case is $N^{-2(1+\beta')/(1+2(1+\beta'))}$, which is smaller than $N^{-2/3}$. We now present a recursive method which attains the optimal rate.

**Remark 2.2** *Bias correction, rate optimality.* The idea here is to achieve a bias correction and the optimal rate by running two ANRE algorithms with different step sizes $\lambda_1$ and $\lambda_2$ in parallel. Let $\hat{\underline{a}}_{t,N}(\lambda_1)$ and $\hat{\underline{a}}_{t,N}(\lambda_2)$ be the ANRE algorithms with step size $\lambda_1$ and $\lambda_2$ respectively, and assume that $\lambda_1 > \lambda_2$. By using (14) for $i = 1, 2$, we have

$$\mathbb{E}\{\hat{\underline{a}}_{t_0,N}(\lambda_i)\} = \underline{a}(u_0) - \frac{1}{N\lambda_i}F(u_0)^{-1}\dot{\underline{a}}(u_0) + O\left(\frac{1}{(N\lambda_i)^{1+\beta'}}\right). \qquad (19)$$

Since $\underline{a}(u_0) - (N\lambda_i)^{-1}F(u_0)^{-1}\dot{\underline{a}}(u_0) \approx \underline{a}(u_0 - (N\lambda_i)^{-1}F(u_0)^{-1})$, we heuristically estimate $\underline{a}(u_0 - (N\lambda_i)^{-1}F(u_0)^{-1})$ instead of $\underline{a}(u_0)$ by the algorithm. By using two different $\lambda_i$ we can find a linear combination of the corresponding estimates such that we 'extrapolate' the two values $\underline{a}(u_0) - (N\lambda_i)^{-1}F(u_0)^{-1}\dot{\underline{a}}(u_0)$ $(i = 1, 2)$ to $\underline{a}(u_0)$. Formally, let $0 < w < 1$, $\lambda_2 = w\lambda_1$ and

$$\tilde{\underline{a}}_{t_0,N}(w) = \frac{1}{1-w}\hat{\underline{a}}_{t_0,N}(\lambda_1) - \frac{w}{1-w}\hat{\underline{a}}_{t_0,N}(\lambda_2).$$

If $|t_0/N - u_0| < 1/N$, then by using (19) we have

$$\mathbb{E}\{\underline{\hat{a}}_{t_0,N}(w)\} = \underline{a}(u_0) + O\left(\frac{1}{(N\lambda)^{1+\beta'}}\right).$$

By using Propostion 4.3 we have

$$\mathbb{E}|\underline{\hat{a}}_{t_0,N} - \underline{a}(u_0)|^2 = O\left(\lambda + \frac{1}{(N\lambda)^{2(1+\beta')}}\right),$$

and choosing $\lambda = \text{const.} \times N^{-(2+2\beta')/(3+2\beta')}$ gives the optimal rate. There remains the problem of choosing $\lambda$ (and $w$). It is obvious that $\lambda$ should be chosen adaptively to the degree of non-stationarity. That is, $\lambda$ should be large if the characteristics of the process are changing more rapidly. However, a more specific suggestion would require more investigations – both theoretically and by simulation.

The above method cannot be extended to higher-order derivatives, since the other remainders are of a lower order $(N\lambda)^{-2}$ (see (15) and the proof of Theorem 2.3).

Finally, we mention that choosing $\lambda_2 < w\lambda_1$ will lead to an estimator of $\underline{a}(u_0 + \Delta)$ with some $\Delta > 0$ (with rate as above). This could be the basis for the prediction of volatility of tvARCH processes. □

# 3. Practical implications

Suppose that we observe data from a (non-rescaled) tvARCH process in discrete time

$$X_t = Z_t\sigma_t, \qquad \sigma_t^2 = a_0(t) + \sum_{j=1}^p a_j(t)X_{t-j}^2, \qquad t \in \mathbb{Z}. \tag{20}$$

In order to estimate $\underline{a}(t)$ we use the estimator $\underline{\hat{a}}_t$ as given in (2) (with all subscripts $N$ dropped). An approximation for the distribution of the estimator is given by Theorem 2.2. Theorem 2.2(ii) can be used directly since it is completely formulated without $N$. The matrices $F(u_0)$ and $\Sigma(u_0)$ depend on the unknown stationary approximation $\mathcal{X}_t(u_0)$ of the process at $u_0 = t_0/N$, that is, at time $t_0$ in non-rescaled time. Since this approximation is unknown we may instead use the process itself in a small neighbourhood of $t_0$, that is, we may estimate, for example, $F(u_0)$ by

$$\frac{1}{m}\sum_{j=0}^{m-1}\frac{\mathcal{X}_{t_0-j}\mathcal{X}_{t_0-j}^\top}{|\mathcal{X}_{t_0-j}|_1^2}$$

with $m$ small and $\mathcal{X}_{t-1}^\top = (1, X_{t-1}^2, \ldots, X_{t-p}^2)$. An estimator which fits the recursive algorithm better is

$$[1 - (1-\lambda)^{t_0-p+1}]^{-1}\sum_{j=0}^{t_0-p}\lambda(1-\lambda)^j\frac{\mathcal{X}_{t_0-j}\mathcal{X}_{t_0-j}^\top}{|\mathcal{X}_{t_0-j}|_1^2}.$$

In the same way we can estimate $\Sigma(u_0)$ which altogether leads, for example, to an approximate confidence interval for $\underline{\hat{a}}_t$. In a similar way Theorem 2.2(i) can be used.

The situation is more difficult with Theorems 2.1 and 2.3, since here the results depend (at first sight) on $N$. Suppose that we have parameter functions $\tilde{a}_j(\cdot)$ and some $N > t_0$ with $\tilde{a}_j(t_0/N) = a_j(t_0)$ (i.e. the original function has been rescaled to the unit interval). Consider Theorem 2.3 with the functions $\tilde{a}_j(\cdot)$. The bias in (14) and (15) contains the term

$$\frac{1}{N}\dot{\tilde{a}}_j(u_0) \approx \frac{1}{N} \frac{\tilde{a}_j(t_0/N) - \tilde{a}_j((t_0-1)/N)}{1/N} = a_j(t_0) - a_j(t_0-1),$$

which again is independent of $N$. To avoid confusion we mention that $N^{-1}\dot{\tilde{a}}_j(u_0)$ of course depends on $N$ once the function $\tilde{a}_j(\cdot)$ has been fixed (as in the asymptotic approach of this paper) but it does not depend on $N$ when it is used to approximate the function $a_j(t)$ since then the function $\tilde{a}_j(\cdot)$ is a different one for each $N$. In the spirit of the remarks above we would, for example, use the expression

$$[1 - (1-\lambda)^{t_0-p+1}]^{-1} \sum_{j=0}^{t_0-p} \lambda(1-\lambda)^j [a_j(t_0) - a_j(t_0-j)]$$

as an estimator of $N^{-1}\dot{\tilde{a}}_j(u_0)$ in (14) and (15).

These considerations also demonstrate the need for the asymptotic approach of this paper. While it is not possible to set down a meaningful asymptotic theory for the model (20) and to derive, for example, a central limit theorem for the estimator $\underline{\hat{a}}_t$, the approach of the present paper for the rescaled model (1) leads to such results. This is achieved by the 'infill asymptotics' where more and more data become available for each local structure (e.g. about time $u_0$) as $N \to \infty$. The results can then be used also for approximations in the model (20) – for example, for confidence intervals.

## 4. Proofs

### 4.1. Some preliminary results

In the next lemma we give a bound for the approximation error between $X_{t,N}^2$ and $X_t(u)^2$. The proofs of these results and further details can be found in Dahlhaus and Subba Rao [4]; see also Subba Rao [14].

**Lemma 4.1.** *Suppose Assumption 2.1 holds with some $r \geq 1$. Let $\{X_{t,N}\}$ and $\{X_t(u)\}$ be defined as in (1) and (5). Then we have:*

(i) *$\{X_t(u)^2\}_t$ is a stationary ergodic process. Furthermore, there exists a stochastic process $\{V_{t,N}\}_t$ and a stationary ergodic process $\{W_t\}_t$ with $\sup_{t,N} \mathbb{E}(|V_{t,N}|^r) < \infty$* *and* $\mathbb{E}(|W_t|^r) < \infty$, *such that*

$$|X_{t,N}^2 - X_t(u)^2| \leq \frac{1}{N^\beta} V_{t,N} + \left|\frac{t}{N} - u\right|^\beta W_t, \tag{21}$$

$$|X_t(u)^2 - X_t(v)^2| \leq |u - v|^\beta W_t. \tag{22}$$

(ii) $\sup_{t,N} \mathbb{E}(X_{t,N}^{2r}) < \infty$ *and* $\sup_u \mathbb{E}(X_t(u)^{2r}) < \infty$.

We now define the derivative process by $\{\dot{X}_t^2(u)\}_t$, which is almost surely the derivative of the squared ARCH process $\{X_t(u)^2\}_t$ (i.e., $\dot{X}_t^2(u) = (\partial X_t(u)^2/\partial u)!$).

**Lemma 4.2.** *Suppose Assumption* 2.2 *holds with some* $r \geq 1$. *Then the derivative process* $\{\dot{X}_t^2(u)\}_t$ *is a well-defined stationary ergodic process which satisfies the representation*

$$\dot{X}_t^2(u) = \left\{ \dot{a}_0(u) + \sum_{j=1}^p [a_j(u)\dot{X}_{t-j}^2(u) + \dot{a}_j(u)X_{t-j}(u)^2] \right\} Z_t^2$$

*and*

$$|\dot{X}_t^2(u) - \dot{X}_t^2(v)| \leq |u - v|^{\beta'} W_t, \tag{23}$$

*where* $W_t$ *is the same as in Lemma* 4.1. *Almost surely all paths of* $\{X_t(u)^2\}_u$ *belong to* Lip(1) *and we have the Taylor series expansion*

$$X_{t,N}^2 = X_t(u)^2 + \left(\frac{t}{N} - u\right)\dot{X}_t^2(u) + \left(\left|\frac{t}{N} - u\right|^{1+\beta'} + \frac{1}{N}\right)R_{t,N},$$

*where* $|R_{t,N}|_1 \leq (V_{t,N} + W_t)$. *Furthermore,* $\sup_N X_{t,N}^2 \leq W_t$, $\sup_u X_t(u)^2 \leq W_t$ *and* $\sup_u |\dot{X}_t^2(u)| < W_t$ *with bounded norms* $\sup_u \mathbb{E}(|\dot{X}_t^2(u)|^r) < \infty$ *and* $\sup_{t,N} \mathbb{E}(|R_{t,N}|^r) < \infty$.

Let $\mathcal{F}_t = \sigma(Z_t^2, Z_{t-1}^2, \ldots)$. We have $\mathcal{F}_t = \sigma(X_{t,N}^2, X_{t-1,N}^2, \ldots) = \sigma(X_t(u)^2, X_{t-1}(u)^2, \ldots)$, since a Volterra expansion gives $X_{t,N}^2$ in terms of $\{Z_t^2\}_t$ and the ARCH equations give $Z_t^2$ in terms of $\{X_{t,N}^2\}_t$. We now consider the mixing properties of functions of the processes $\{\mathcal{X}_{t,N}\}_t$ and $\{\mathcal{X}_t(u)\}_t$. The proof of the proposition below can be found in Section A.1.

**Proposition 4.1.** *Suppose Assumption* 2.1 *holds with* $r = 1$. *Let* $\{X_{t,N}\}$ *and* $\{X_t(u)\}$ *be defined as in* (1) *and* (5), *respectively. Then there exists a* $(1 - \eta) \leq \rho < 1$ *such that for any* $\phi \in$ Lip(1),

$$|\mathbb{E}[\phi(\mathcal{X}_{t,N})|\mathcal{F}_{t-k}] - \mathbb{E}[\phi(\mathcal{X}_{t,N})]|_1 \leq K\rho^k(1 + |\mathcal{X}_{t-k,N}|_1), \tag{24}$$

$$|\mathbb{E}[\phi(\mathcal{X}_t(u))|\mathcal{F}_{t-k}] - \mathbb{E}[\phi(\mathcal{X}_t(u))]|_1 \leq K\rho^k(1 + |\mathcal{X}_{t-k}(u)|_1), \tag{25}$$

$$|\mathbb{E}[\phi(\mathcal{X}_{t,N})|\mathcal{F}_{t-k}]|_1 \leq K(1 + \rho^k|\mathcal{X}_{t-k,N}|_1); \tag{26}$$

*if* $j, k > 0$ *then*

$$|\mathbb{E}[\phi(\mathcal{X}_{t,N})|\mathcal{F}_{t-k}] - \mathbb{E}[\phi(\mathcal{X}_{t,N})|\mathcal{F}_{t-k-j}]|_1 \leq K\rho^k(|\mathcal{X}_{t-k,N}|_1 + |\mathcal{X}_{t-k-j,N}|_1); \quad (27)$$

*and if Assumption 2.1 holds with some* $r > 1$ *then for* $1 \leq q \leq r$ *we have*

$$\mathbb{E}(|\mathcal{X}_{t,N}|_1^q|\mathcal{F}_{t-k}) \leq K|\mathcal{X}_{t-k,N}|_1^q, \quad (28)$$

*where the constant* $K$ *is independent of* $t, k, j$ *and* $N$.

The corollary below follows by using (24) and (25).

**Corollary 4.1.** *Suppose Assumption 2.1 holds. Let* $\{X_{t,N}\}$ *and* $\{X_t(u)\}$ *be defined as in* (1) *and* (5) *respectively, and* $\phi \in \mathrm{Lip}(1)$. *Then* $\{\phi(\mathcal{X}_{t,N})\}$ *and* $\{\phi(\mathcal{X}_t(u))\}$ *are* $L_q$-*mixingales of size* $-\infty$.

## 4.2. The pertubation technique

In this section we use the pertubation technique, introduced in Aguech *et al.* [1], to show consistency of the ANRE estimator. To analyse the algorithm we compare it with a similar one driven by the true parameters where $X_{t,N}$ has been replaced by the stationary process $X_t(u)$. Let $\delta_{t,N}(u) = \hat{\underline{a}}_{t,N} - \underline{a}(u)$,

$$\mathcal{M}_{t,N} = (Z_t^2 - 1)\sigma_{t,N}^2 \frac{\mathcal{X}_{t-1,N}}{|\mathcal{X}_{t-1,N}|_1^2}, \qquad \mathcal{M}_t(u) = (Z_t^2 - 1)\sigma_t(u)^2 \frac{\mathcal{X}_{t-1}(u)}{|\mathcal{X}_{t-1}(u)|_1^2}, \quad (29)$$

$$F_{t,N} = \frac{\mathcal{X}_{t,N}\mathcal{X}_{t,N}^\top}{|\mathcal{X}_{t,N}|_1^2}, \qquad\qquad F_t(u) = \frac{\mathcal{X}_t(u)\mathcal{X}_t(u)^\top}{|\mathcal{X}_t(u)|_1^2}, \quad (30)$$

and $F(u)$ be defined as in (9). The 'true algorithm' is

$$\underline{a}(u) = \underline{a}(u) + \lambda\{X_t(u)^2 - \underline{a}(u)^\top \mathcal{X}_{t-1}(u)\}\frac{\mathcal{X}_{t-1}(u)}{|\mathcal{X}_{t-1}(u)|_1^2} - \lambda\mathcal{M}_t(u). \quad (31)$$

An advantage of the specific form of the random matrices $F_{t,N}$ is that $|F_{t,N}|_1 \leq (p+1)^2$. This upper bound will make some of the analysis easier to handle.

By subtracting (31) from (2) we obtain

$$\delta_{t,N}(u) = (I - \lambda F_{t-1,N})\delta_{t-1,N}(u) + \lambda\mathcal{B}_{t,N}(u) + \lambda\mathcal{M}_{t,N}(u), \quad (32)$$

where

$$\mathcal{B}_{t,N}(u) = \{\mathcal{M}_{t,N} - \mathcal{M}_t(u)\} + F_{t-1,N}\left\{\underline{a}\left(\frac{t}{N}\right) - \underline{a}(u)\right\}. \quad (33)$$

We note that (32) can also be considered as a recursive algorithm, thus we call (32) the error algorithm. There are two terms driving the error algorithm: the bias $\mathcal{B}_{t,N}(u)$ and the stochastic term $\mathcal{M}_t(u)$. Because the error algorithm is linear with respect to the estimators we can separate $\delta_{t,N}(u)$ in terms of the bias and the stochastic terms:

$$\delta_{t,N}(u) = \delta_{t,N}^B(u) + \delta_{t,N}^M(u) + \delta_{t,N}^R(u),$$

where

$$\delta_{t,N}^B(u) = (I - \lambda F_{t-1,N})\delta_{t-1,N}^B(u) + \lambda \mathcal{B}_{t,N}(u), \qquad \delta_{p,N}^B(u) = 0, \tag{34}$$

$$\delta_{t,N}^M(u) = (I - \lambda F_{t-1,N})\delta_{t-1,N}^M(u) + \lambda \mathcal{M}_t(u), \qquad \delta_{p,N}^M(u) = 0, \tag{35}$$

$$\delta_{t,N}^R(u) = (I - \lambda F_{t-1,N})\delta_{t-1,N}^R(u), \qquad\qquad \delta_{p,N}^R(u) = -\underline{a}(u). \tag{36}$$

We have for $y \in \{B, M, R\}$,

$$\delta_{t,N}^y(u) = \sum_{k=0}^{t-p-1} \left\{ \prod_{i=0}^{k-1} (I - \lambda F_{t-i-1,N}) \right\} \mathcal{D}_{t-k,N}^y(u), \tag{37}$$

where $\mathcal{D}_{t,N}^B(u) = \lambda \mathcal{B}_{t,N}(u)$, $\mathcal{D}_{t,N}^M(u) = \lambda \mathcal{M}_t(u)$, $\mathcal{D}_{t,N}^R(u) = 0$ if $t > p$ and $\mathcal{D}_{p,N}^R(u) = -\underline{a}(u)$.

For $\delta_{t,N}^y(u)$ to converge, it is clear that the random product $\prod_{i=0}^{t-1}(I - \lambda F_{t-i-1,N})$ must decay. Technically this is one of the main results. It is stated in the following theorem.

**Theorem 4.1.** *Suppose Assumption* 2.1 *holds with* $r = 4$ *and* $N$ *is sufficiently large. Then for all* $q \geq 1$ *and* $(p+1) \leq t \leq N$ *there exists* $M > 0$, $\delta > 0$ *such that*

$$\mathbb{E} \left\| \prod_{i=0}^{k-1} \left( I - \lambda \frac{\mathcal{X}_{t-i-1,N} \mathcal{X}_{t-i-1,N}^\top}{|\mathcal{X}_{t-i-1,N}|_1^2} \right) \right\|^q \leq M \exp\{-\delta \lambda k\}. \tag{38}$$

Under Assumption 2.1 and by using the proposition above, there exists a $\delta > 0$ with

$$(\mathbb{E}|\delta_{t,N}^R|^q)^{1/q} \leq M \exp\{-\delta \lambda t\}, \tag{39}$$

for $1 \leq q \leq r$ and $t = p+1, \ldots, N$. Therefore this term decays exponentially and, as we shall see below, is of lower order than $\delta_{t,N}^B(u)$ and $\delta_{t,N}^M(u)$.

We now study the stochastic bias at $t_0$ with $|t_0/N - u_0| < 1/N$. It is straightforward to see that $\delta_{t,N}^B(u) = \delta_{t,N}^{1,B}(u) + \delta_{t,N}^{2,B}(u) + \delta_{t,N}^{3,B}(u)$, where

$$\delta_{t,N}^{1,B}(u) = (I - \lambda F_{t-1,N})\delta_{t-1,N}^{1,B}(u) + \lambda \{\mathcal{M}_{t,N} - \mathcal{M}_t(u)\},$$

$$\delta_{t,N}^{2,B}(u) = (I - \lambda F_{t-1,N})\delta_{t-1,N}^{2,B}(u) + \lambda F(u) \left\{ \underline{a}\left(\frac{t}{N}\right) - \underline{a}(u) \right\}, \tag{40}$$

$$\delta_{t,N}^{3,B}(u) = (I - \lambda F_{t-1,N})\delta_{t-1,N}^{3,B}(u) + \lambda (F_{t-1,N} - F(u)) \left\{ \underline{a}\left(\frac{t}{N}\right) - \underline{a}(u) \right\},$$

$\delta_{p,N}^{1,B} = 0$, $\delta_{p,N}^{2,B} = 0$ and $\delta_{p,N}^{3,B} = 0$. In order to obtain asymptotic expressions for the expectation of each component $\delta_{t_0,N}^{1,B}$, $\delta_{t_0,N}^{2,B}$, $\delta_{t_0,N}^{2,B}$ and $\delta_{t_0,N}^{M}$ we use the pertubation technique proposed by Aguech *et al.* [1]. For $x = M, (1, B), (2, B)$, we can decompose the stochastic and bias terms as follows:

$$\delta_{t_0,N}^x(u_0) = J_{t_0,N}^{x,1}(u_0) + J_{t_0,N}^{x,2}(u_0) + H_{t_0,N}^x(u_0), \tag{41}$$

where

$$
\begin{aligned}
J_{t,N}^{x,1}(u_0) &= (I - \lambda F(u_0))J_{t-1,N}^{x,1}(u_0) + G_{t,N}^x, \\
J_{t,N}^{x,2}(u_0) &= (I - \lambda F(u_0))J_{t-1,N}^{x,2}(u_0) - \lambda\{F_{t-1,N} - F(u_0)\}J_{t-1,N}^{x,1}(u_0), \\
H_{t,N}^x(u_0) &= (I - \lambda F_{t-1,N})H_{t-1,N}^x(u_0) - \lambda\{F_{t-1,N} - F(u_0)\}J_{t-1,N}^{x,2}(u_0),
\end{aligned}
\tag{42}
$$

with, for $t < t_0$,

$$
\begin{aligned}
G_{t,N}^M &= \lambda\mathcal{M}_t(u_0), & G_{t,N}^{1,B} &= \lambda[\mathcal{M}_{t,N} - \mathcal{M}_t(u_0)], \\
G_{t,N}^{2,B} &= F(u_0)\Big[\underline{a}\Big(\frac{t}{N}\Big) - \underline{a}(u)\Big], & G_{t,N}^{3,B} &= (F_{t-1,N} - F(u_0))\Big[\underline{a}\Big(\frac{t}{N}\Big) - \underline{a}(u)\Big].
\end{aligned}
\tag{43}
$$

Furthermore, $J_{p,N}^{x,1}(u_0) = J_{p,N}^{x,2}(u_0) = H_{p,N}^x(u_0) = 0$. Equation (41) can easily be derived by taking the sum on both sides of the three equations in (42). In the proposition below we will show that $J_{t_0,N}^{x,1}(u_0)$ is the principal term in the expansion of $\delta_{t_0,N}^x(u_0)$, that is, with $\delta_{t_0,N}^x(u_0) \approx J_{t_0,N}^{x,1}(u_0)$. Substituting (43) into (42) gives

$$
\begin{aligned}
J_{t_0,N}^{(1,B),1}(u_0) &= \sum_{k=0}^{t_0-p-1} \lambda\{I - \lambda F(u_0)\}^k\{\mathcal{M}_{t_0-k,N} - \mathcal{M}_{t_0-k}(u_0)\}, \\
J_{t_0,N}^{(2,B),1}(u_0) &= \sum_{k=1}^{t_0-p-1} \{I - \lambda F(u_0)\}^k F(u_0)\Big\{\underline{a}\Big(\frac{t_0-k}{N}\Big) - \underline{a}(u_0)\Big\}, \\
J_{t_0,N}^{(3,B),1}(u_0) &= \sum_{k=1}^{t_0-p-1} \{I - \lambda F(u_0)\}^k(F_{t-k-1,N} - F(u_0))\Big\{\underline{a}\Big(\frac{t_0-k}{N}\Big) - \underline{a}(u_0)\Big\}.
\end{aligned}
\tag{44}
$$

In the proof below we require the following definition

$$D_{t,N} = \sum_{i=0}^{p-1}[V_{t-i,N} + W_{t-i}], \qquad \text{for } t \geq p, \tag{45}$$

where $V_{t,N}$ and $W_t$ are defined as in Lemma 4.1.

**Proposition 4.2.** *Suppose Assumption 2.1 holds with $r > 4$, and let $\delta_{t_0,N}^B(u_0)$, $J_{t_0,N}^{(1,B),1}(u_0)$, $J_{t_0,N}^{(2,B),1}(u_0)$, $J_{t_0,N}^{(3,B),1}(u_0)$, be defined as in (34) and (44). Then for*

$|t_0/N - u_0| < 1/N$ and $\lambda N \geq (\log N)^{1+\varepsilon}$, for some $\varepsilon > 0$, we have

(i)  $(\mathbb{E}|J_{t_0,N}^{(1,B),1}(u_0) + J_{t_0,N}^{(2,B),1}(u_0)|^r)^{1/r} = O\left(\dfrac{1}{(N\lambda)^\beta}\right),$  (46)

(ii)  $(\mathbb{E}|J_{t_0,N}^{(3,B),1}(u_0)|^r)^{1/r} = O\left(\dfrac{1}{(N\lambda)^{2\beta}} + \dfrac{\lambda^{1/2}}{(N\lambda)^\beta}\right),$

(iii)  $(\mathbb{E}|\delta_{t_0,N}^B(u_0) - J_{t_0,N}^{(1,B),1}(u_0) - J_{t_0,N}^{(2,B),1}(u_0)|^2)^{1/2} = O\left(\dfrac{1}{(N\lambda)^{2\beta}} + \dfrac{\sqrt{\lambda}}{(N\lambda)^\beta} + \lambda\right).$  (47)

**Proof.** We first prove (i). Let us consider $J_{t_0,N}^{(1,B),1}(u_0)$. By using (105) we have

$$|\mathcal{M}_{t,N} - \mathcal{M}_t(u_0)|_1 \leq |Z_t^2 + 1|\left\{\frac{1}{N^\beta} + \left(\frac{p}{N}\right)^\beta + \left|\frac{t}{N} - u_0\right|^\beta\right\}D_{t,N}.$$

Therefore, by substituting the above bound into $J_{t_0,N}^{(1,B),1}(u_0)$ and by using (99), we have

$$(\mathbb{E}|J_{t_0,N}^{(1,B),1}(u_0)|^r)^{1/r}$$
$$\leq \frac{K}{N^\beta}\sum_{k=0}^{t_0-p-1}\lambda(1-\delta\lambda)^k\{1+p^\beta+k^\beta\}(\mathbb{E}(Z_0^2+1)^r)^{1/r}(\mathbb{E}|D_{t_0-k-1,N}|^r)^{1/r}.$$

Now by using Lemma 4.1(i) we have that $\sup_{t,N}|D_{t,N}|^r < \infty$. Furthermore, from Lemma C.3 in Moulines *et al.* [11] we have

$$\sum_{k=1}^{N}(1-\lambda)^k k^\beta \leq \lambda^{-1-\beta}.$$  (48)

By using the above we obtain

$$(\mathbb{E}|J_{t_0,N}^{(1,B),1}(u_0)|^r)^{1/r} \leq K\sup_{t,N}(\mathbb{E}|D_{t,N}|^r)^{1/r}(\mathbb{E}|Z_0^2+1|^r)^{1/r}\frac{1}{(N\lambda)^\beta} = O\left(\frac{1}{(N\lambda)^\beta}\right).$$  (49)

We now bound $(\mathbb{E}|J_{t_0,N}^{(2,B),1}(u_0)|^r)^{1/r}$. Since $\underline{a}(u) \in \text{Lip}(\beta)$, by using (99) and (48) we have

$$(\mathbb{E}|J_{t_0,N}^{(2,B),1}(u_0)|^r)^{1/r} \leq C\sum_{k=1}^{t_0-p-1}\lambda(1-\delta\lambda)^k\left(\frac{k}{N}\right)^\beta = O\left(\frac{1}{(\lambda N)^\beta}\right).$$  (50)

Thus (49) and (50) give the bound (46), which completes the proof of (i).

To prove (ii), we now bound $(\mathbb{E}|J_{t_0,N}^{(3,B),1}(u_0)|^r)^{1/r}$. We first observe that $J_{t_0,N}^{(3,B),1}(u_0)$ can be written as

$$J_{t_0,N}^{(3,B),1}(u_0) = I_{t_0,N} + II_{t_0,N},$$

where

$$I_{t_0,N} = \sum_{k=1}^{t_0-p-1} \{I - \lambda F(u_0)\}^k (F_{t-k-1,N} - F_{t-k-1}(u_0)) \left\{ \underline{a}\left(\frac{t_0-k}{N}\right) - \underline{a}(u_0) \right\}$$

$$II_{t_0,N} = \sum_{k=1}^{t_0-p-1} \{I - \lambda F(u_0)\}^k (F_{t-k-1}(u_0) - F(u_0)) \left\{ \underline{a}\left(\frac{t_0-k}{N}\right) - \underline{a}(u_0) \right\}.$$

Using the same proof as for $J_{t_0,N}^{(1,B),1}(u_0)$, it is straightforward to show that $(\mathbb{E}|I_{t_0,N}|^r)^{1/r} = O((N\lambda)^{-2\beta})$. In order to bound $II_{t_0,N}$, let $\bar{F}_k(u_0) = F_k(u_0) - F(u_0)$ and $\phi(\underline{x}) = \underline{x}\underline{x}^\top/|\underline{x}|_1^2$, where $\underline{x} = (1, x_1, \ldots, x_p)$; then $\phi(\underline{x}) \in \mathrm{Lip}(1)$ and $F_{t,N} = \phi(\mathcal{X}_{t,N})$. Since $\phi \in \mathrm{Lip}(1)$, by using Corollary 4.1 we have that $F_t(u_0) - F(u_0)$ can be written as the sum of martingale differences $F_t(u_0) - F(u_0) = \sum_{\ell=0}^{\infty} m_t(\ell)$, where $m_t(\ell)$ is a $(p+1) \times (p+1)$ matrix defined by $m_t(\ell) = \{\mathbb{E}(F_{t,N}|\mathcal{F}_{t-\ell}) - \mathbb{E}(F_{t,N}|\mathcal{F}_{t-\ell-1})\}$. By using (27) we have $(\mathbb{E}|m_t(\ell)|^{r/2})^{2/r} \le K\rho^\ell$. Substituting the above into $B_{t_0,N}^{2,B}$ we have

$$II_{t_0,N} = \sum_{\ell=0}^{\infty} \sum_{k=0}^{t_0-p-1} \lambda \{I - \lambda F(u_0)\}^k m_{t_0-k-1}(\ell) \left\{ \underline{a}\left(\frac{t_0-k}{N}\right) - \underline{a}(u_0) \right\}.$$

We note that $\|\{I - \lambda F(u_0)\}^k\| \le K(1-\lambda\delta)^k$ (see (99)). Furthermore, if $t_1 < t_2$ then $\mathbb{E}\{m_{t_1}(\ell)m_{t_2}(\ell)\} = \mathbb{E}\{m_{t_1}(\ell)\mathbb{E}(m_{t_2}(\ell)|\mathcal{F}_{t_1-\ell})\} = 0$, therefore $\{m_t(\ell)\}_t$ is a sequence of martingale differences. Since $J_{t_0-k-1,N}^{(2,B),1}$ is deterministic we can use Burkhölder's inequality (cf. Hall and Heyde [6], Theorem 12.2) to obtain

$$(\mathbb{E}|II_{t_0,N}|^r)^{1/r} \le K\lambda \sum_{\ell=0}^{\infty} \left( \mathbb{E} \left| \sum_{k=0}^{t_0-p-1} \{I - \lambda F(u_0)\}^k m_{t_0-k-1}(\ell) \left\{ \underline{a}\left(\frac{t_0-k}{N}\right)_j - \underline{a}(u_0)_j \right\} \right|^r \right)^{1/r}$$

$$\le K\lambda \sum_{\ell=0}^{\infty} \left\{ 2r \sum_{k=0}^{t_0-p-1} (1-\lambda\delta)^{2k} (\mathbb{E}|\{m_k(\ell)\}_{ij}|^r)^{2/r} \left(\frac{k}{N}\right)^\beta \right\}$$

$$\le \frac{K}{(N\lambda)^\beta} \sum_{\ell=0}^{\infty} \rho^\ell \lambda^{1/2} = O\left(\frac{\lambda^{1/2}}{(N\lambda)^\beta}\right). \tag{51}$$

Thus we have proved (ii).

   We now prove (iii). By using (41) we have

$$(\mathbb{E}|\delta_{t_0,N}^B(u_0) - J_{t_0,N}^{(1,B),1}(u_0) - J_{t_0,N}^{(2,B),1}(u_0)|^{r/2})^{2/r}$$

$$= (\mathbb{E}|J_{t_0,N}^{(3,B),1}(u_0)|^{r/2})^{2/r} + \left( \mathbb{E} \left| \sum_{i=1}^{3} J_{t_0,N}^{(i,B),2}(u_0) \right|^{r/2} \right)^{2/r}$$

$$+ \left( \mathbb{E} \left| \sum_{i=1}^{3} H_{t_0,N}^{(i,B)}(u_0) \right|^{r/2} \right)^{2/r} + O\left( \frac{1}{N^\beta} + \exp\{-\lambda \delta t_0\} \right).$$

We have bounded the first term of the above; to bound the rest we partition $\sum_{i=1}^{3} J_{t_0,N}^{(i,B),2}$ into four terms:

$$\sum_{i=1}^{3} J_{t_0,N}^{(i,B),2} = A_{t_0,N}^{B} + B_{t_0,N}^{1,B} + B_{t_0,N}^{2,B} + B_{t_0,N}^{3,B}, \tag{52}$$

where

$$A_{t_0,N}^{B} = - \sum_{k=0}^{t_0-p-1} \lambda \{I - \lambda F(u_0)\}^k [F_{t_0-k-1,N} - F_{t_0-k-1}(u_0)] \times \left( \sum_{i=1}^{3} J_{t_0-k-1,N}^{(i,B),1} \right)$$

and, for $y \in \{1,2,3\}$,

$$B_{t_0,N}^{y,B} = - \sum_{k=0}^{t_0-p-1} \lambda \{I - \lambda F(u_0)\}^k [F_{t_0-k-1}(u_0) - F(u_0)] J_{t_0-k-1,N}^{(y,B),1}(u_0).$$

We first bound $A_{t_0,N}^{B}$. By using (49)–(51) and (107) we have

$$(\mathbb{E}|A_{t_0,N}^{B}|^{r/2})^{2/r} \leq \frac{K}{(N\lambda)^{2\beta}}.$$

By using (109) and (110) we have

$$(\mathbb{E}|B_{t_0,N}^{1,B}|^{r/2})^{2/r} = \left( \mathbb{E} \left| \sum_{k=0}^{t_0-p-1} \sum_{i=0}^{t_0-p-k-1} \lambda^2 \{I - \lambda F(u_0)\}^{k+i} \bar{F}_{t_0-k-1}(u_0) \right. \right.$$
$$\left. \left. \times \{\mathcal{M}_{t_0-k-i-1,N} - \mathcal{M}_{t_0-k-i-1}(u_0)\} \right|^{r/2} \right)^{2/r}$$

$$\leq K\lambda.$$

Using a similar proof to the above to bound $(\mathbb{E}|J_{t_0-k-1,N}^{(3,B),1}(u_0)|^r)^{1/r}$ we can show that $\|B_{t_0,N}^{2,B}\|_{r/2}^{E} = O(\lambda^{1/2}(N\lambda)^{-\beta})$. By using (51) it is straightforward to show that $\|B_{t_0,N}^{3,B}\|_{r/2}^{E} = O((N\lambda)^{-2\beta} + \lambda^{1/2}(N\lambda)^{-\beta})$. Substituting the above bounds for $(\mathbb{E}|A_{t_0,N}^{B}|^{r/2})^{2/r}$, $(\mathbb{E}|B_{t_0,N}^{1,B}|^{r/2})^{2/r}$, $(\mathbb{E}|B_{t_0,N}^{2,B}|^{r/2})^{2/r}$ and $(\mathbb{E}|B_{t_0,N}^{3,B}|^{r/2})^{2/r}$ into (52), we obtain

$$(\mathbb{E}|J_{t_0,N}^{(1,B),2}(u_0) + J_{t_0,N}^{(2,B),2}(u_0) + J_{t_0,N}^{(3,B),2}(u_0)|^{r/2})^{2/r} = O(\delta_N). \tag{53}$$

Finally, we prove (iii) by bounding, for $y = 1,2,3$, $\|H_{t_0,N}^{(y,B)}(u_0)\|_2^{E}$. By using Hölder's inequality, (53), Theorem 4.1 and that $\{F_{t,N}\}$ are bounded random matrices, we have

$$
\begin{aligned}
&(\mathbb{E}|H_{t_0,N}^{(1,B)}(u_0) + H_{t_0,N}^{(2,B)}(u_0) + H_{t_0,N}^{(3,B)}(u_0)|^2)^{1/2} \\
&\leq \sum_{k=0}^{t_0-p-1} \lambda \left( \mathbb{E} \left\| \left( \prod_{i=0}^{k-1} (I - \lambda F_{t_0-i,N}) \right) \right\|^{2r/(r-4)} \right)^{(r-4)/2r} \\
&\qquad \times \left( \mathbb{E}|[F_{t_0-k-1,N} - F(u_0)]\{J_{t_0-k-1,N}^{(1,B),2}(u_0) + J_{t_0-k-1,N}^{(2,B),2}(u_0) + J_{t_0-k-1,N}^{(3,B),2}\}|^{r/2} \right)^{2/r} \\
&\leq K \left( \frac{1}{(N\lambda)^{2\beta}} + \frac{\lambda^{1/2}}{(N\lambda)^{\beta}} + \lambda \right).
\end{aligned}
\tag{54}
$$

Since $(\mathbb{E}|J_{t_0,N}^{(i,B),2}(u_0)|^2)^{1/2} \leq (\mathbb{E}|J_{t_0,N}^{(i,B),2}(u_0)|^{r/2})^{2/r}$, by substituting (53) and (54) into (52) we have

$$
(\mathbb{E}|\delta_{t_0,N}^{B}(u_0) - J_{t_0,N}^{(1,B),1}(u_0) - J_{t_0,N}^{(2,B),1}(u_0)|^2)^{1/2} = O\left( \frac{1}{(N\lambda)^{2\beta}} + \frac{\sqrt{\lambda}}{(N\lambda)^{\beta}} + \lambda \right).
$$

and we obtain (47). □

In the following lemma we show that $\delta_{t_0,N}^{M}(u_0)$ is dominated by $J_{t_0,N}^{M,1}(u_0)$, where

$$
J_{t_0,N}^{M,1}(u_0) = \sum_{k=0}^{t_0-p-1} \lambda (I - \lambda F(u_0))^k \mathcal{M}_{t_0-k}(u_0).
\tag{55}
$$

We observe that $J_{t_0,N}^{M,1}(u_0)$ and $\mathcal{L}_{t_0}(u_0)$ (defined in (8)) are the same.

**Proposition 4.3.** *Suppose Assumption 2.1 holds with $r > 4$. Let $\delta_{t_0,N}^{M}(u_0)$ and $J_{t_0,N}^{M,1}(u_0)$ be defined as in (35) and (55). Then for $|t_0/N - u_0| < 1/N$ and $\lambda N \geq (\log N)^{1+\varepsilon}$, for some $\varepsilon > 0$, we have*

$$
(\mathbb{E}|J_{t_0,N}^{M,1}(u_0)|^r)^{1/r} = O(\sqrt{\lambda}),
\tag{56}
$$

$$
(\mathbb{E}|\delta_{t_0,N}^{M}(u_0) - J_{t_0,N}^{M,1}(u_0)|^2)^{1/2} = O\left( \lambda + \frac{\sqrt{\lambda}}{(N\lambda)^{\beta}} \right).
\tag{57}
$$

**Proof.** Since each component of the vector sequence $\{\mathcal{M}_t(u_0)\}$ is a martingale difference, we can use the Burkhölder inequality componentwise and (99) to obtain (56). Since $\delta_{t_0,N}^{M}(u_0) - J_{t_0,N}^{M,1}(u_0) = J_{t_0,N}^{M,2}(u_0) + H_{t_0,N}^{M}(u_0)$, to prove (57) we bound $J_{t_0,N}^{M,2}(u_0)$ and $H_{t_0,N}^{M}(u_0)$. By using the same arguments given in Proposition 4.2 we also have

$$
(\mathbb{E}|J_{t_0,N}^{M,2}(u_0)|^{r/2})^{2/r} \leq K \left( \lambda + \frac{\sqrt{\lambda}}{(N\lambda)^{\beta}} \right).
$$

Finally, we obtain a bound for $H_{t_0,N}^M(u_0)$. By using Hölder's inequality, Theorem 4.1 and $(\mathbb{E}\|F_{t_0-k-1,N} - F(u_0)\|^{r/2})^{2/r} \leq 2(p+1)^2$ we have

$$(\mathbb{E}|H_{t_0,N}^M(u_0)|^2)^{1/2} \leq K\left(\lambda + \frac{\sqrt{\lambda}}{(N\lambda)^\beta}\right).$$

Thus we have shown (57). $\qquad\square$

By using Propositions 4.2 and 4.3 we have established that $J_{t_0,N}^{(1,B),1}(u_0)$, $J_{t_0,N}^{(2,B),1}(u_0)$ and $J_{t_0,N}^{M,1}(u_0)$ are the principle terms in $\hat{\underline{a}}_{t_0,N} - \underline{a}(u_0)$. More precisely, we have shown that

$$\hat{\underline{a}}_{t_0,N} - \underline{a}(u_0) = J_{t_0,N}^{(1,B),1}(u_0) + J_{t_0,N}^{(2,B),1}(u_0) + J_{t_0,N}^{M,1}(u_0) + R_N^{(1)}, \tag{58}$$

where $(E|R_N^{(1)}|^2)^{1/2} = O(\delta_N)$ ($\delta_N$ is defined in (7)). By using (58), we show below that an upper bound for the mean squared error can immediately be obtained.

**Proof of Theorem 2.1.** By substituting the bounds for $(\mathbb{E}|J_{t_0,N}^{(1,B),1}(u_0) + J_{t_0,N}^{(2,B),1}(u_0)|^2)^{1/2}$ and $(\mathbb{E}|J_{t_0,N}^{M,1}(u_0)|^2)^{1/2}$ (in Propositions 4.2 and 4.3) into (58) we have $\mathbb{E}(|\hat{\underline{a}}_{t_0,N} - \underline{a}(u_0)|^2) = O((\lambda N)^{-\beta} + \lambda^{1/2} + \delta_N)^2$. Thus we have (10). $\qquad\square$

## 4.3. Proof of Theorem 2.2

We can see from Propositions 4.2 and 4.3 that $J_{t_0,N}^{(1,B),1}(u_0) + J_{t_0,N}^{(1,B),2}(u_0)$ and $J_{t_0,N}^{M,1}(u_0)$ are leading terms in the expansion of $\delta_{t_0,N}$. For this reason, to prove Theorem 2.2 we need only consider these terms. Now considering the bias, we observe that

$$J_{t_0,N}^{(1,B),1}(u_0) + J_{t_0,N}^{(2,B),1}(u_0) = \mathcal{R}_{t_0,N}(u_0) + R_N^{(2)}, \tag{59}$$

where have replaced $\mathcal{M}_{t,N}$ by $\mathcal{M}_t(t/N)$, leading to the remainder

$$R_N^{(2)} = \sum_{k=0}^{t_0-p-1} \lambda\{I - \lambda F(u_0)\}^k \left\{\mathcal{M}_{t_0-k,N} - \mathcal{M}_{t_0-k}\left(\frac{t_0-k}{N}\right)\right\}.$$

By using (21) we have $(\mathbb{E}|R_N^{(2)}|^r)^{1/r} \leq K/N^\beta$, for all $t, N$. Therefore under Assumption 2.1 with $r > 2$, if $\lambda N \geq (\log N)^{1+\varepsilon}$ for some $\varepsilon > 0$, we have, by substituting (59) into (58),

$$\hat{\underline{a}}_{t_0,N} - \underline{a}(u_0) = \mathcal{L}_{t_0}(u_0) + \mathcal{R}_{t_0,N}(u_0) + R_N^{(3)}, \tag{60}$$

where $(\mathbb{E}|R_N^{(3)}|^2)^{1/2} = O(\delta_N)$. In the proposition below we show the asymptotic normality of $\mathcal{L}_{t_0}(u_0)$, which we use together with (60) to obtain the asymptotic normality of $\hat{\underline{a}}_{t_0,N}$.

**Proposition 4.4.** *Suppose Assumption* 2.1 *holds with* $r = 1$ *and, for some* $\delta > 0$, $\mathbb{E}(Z_0^{4+\delta}) < \infty$. *Let* $\mathcal{L}_{t_0}(u_0)$ *and* $\Sigma(u)$ *be defined as in* (8) *and* (13), *respectively. If* $|t_0/N - u_0| < 1/N$, *then we have*

$$\lambda^{-1/2}\mathcal{L}_{t_0}(u_0) \xrightarrow{\mathcal{D}} \mathcal{N}(0, \Sigma(u_0)), \tag{61}$$

*where* $\lambda \to 0$ *as* $N \to \infty$ *and* $\lambda N \geq (\log N)^{1+\varepsilon}$, *for some* $\varepsilon > 0$.

**Proof.** Since $\mathcal{L}_{t_0}(u_0)$ is the weighted sum of martingale differences the result follows from the martingale central limit theorem and the Cramér–Wold device. It is straightforward to check the conditional Lindeberg condition, so we omit the details. We simply note that by using (100) we obtain the limit of the conditional variance of $\mathcal{L}_{t_0}(u_0)$:

$$\mu_4 \sum_{k=0}^{t_0-p-1} \lambda^2 \{I - \lambda F(u_0)\}^{2k} \sigma_{t_0-k}(u_0)^4 \frac{\mathcal{X}_{t_0-k-1}(u_0)\mathcal{X}_{t_0-k-1}(u_0)^\top}{|\mathcal{X}_{t_0-k-1}(u_0)|_1^4} \xrightarrow{\mathcal{P}} \Sigma(u_0). \tag{62}$$

$\square$

**Proof of Theorem 2.2.** It follows from (60) that

$$\lambda^{-1/2}\{\underline{\hat{a}}_{t_0,N} - \underline{a}(u_0)\} = \lambda^{-1/2}\mathcal{R}_{t_0,N}(u_0) + \lambda^{-1/2}\mathcal{L}_{t_0}(u_0) + O_p(\lambda^{1/2}\delta_N).$$

Therefore, if $\lambda \gg N^{-4\beta/(4\beta+1)}$ and $\lambda \gg N^{-2\beta}$ then $\lambda^{-1/2}/(N\lambda)^{2\beta} \to 0$ and $\lambda^{-1/2}/N^\beta \to 0$ respectively, and we have

$$\lambda^{-1/2}\{\underline{\hat{a}}_{t_0,N} - \underline{a}(u_0)\} = \lambda^{-1/2}\mathcal{R}_{t_0,N}(u_0) + \lambda^{-1/2}\mathcal{L}_{t_0}(u_0) + o_p(1).$$

By using the above and Proposition 4.4 we obtain (11). Finally, since $|\mathcal{R}_{t_0,N}(u_0)|_1 = O_p((N\lambda)^{-\beta})$, if $\lambda \gg N^{-2\beta/(2\beta+1)}$ then $\lambda^{-1/2}\mathcal{R}_{t_0,N}(u_0) \xrightarrow{\mathcal{P}} 0$ and we have (12). $\square$

# 5. An asymptotic expression for the bias under additional regularity

Under additional assumptions on the smoothness of the parameters $\underline{a}(u)$, we obtain in this section an exact expression for the bias, which we then use to prove Theorem 2.3.

In the section above we showed that $\delta_{t_0,N}^B \approx \mathcal{R}_{t_0,N}(u_0)$. Since $\mathcal{M}_t(u)$ is a function on $X_t(u)^2$ whose derivative exists, the derivative of $\mathcal{M}_t(u)$ also exists and is given by

$$\dot{\mathcal{M}}_t(u) = (Z_t^2 - 1)\{F_{t-1}(u)\underline{\dot{a}}(u) + \dot{F}_{t-1}(u)\underline{a}(u)\}, \tag{63}$$

where

$$\dot{F}_{t-1}(u) = \underline{Y}_{t-1}(u)\underline{\dot{Y}}_{t-1}(u)^\top + \underline{\dot{Y}}_{t-1}(u)\underline{Y}_{t-1}(u)^\top, \tag{64}$$

with

$$\underline{Y}_{t-1}(u) = \frac{1}{|\mathcal{X}_{t-1}(u)|_1} \mathcal{X}_{t-1}(u),$$

$$\underline{\dot{Y}}_{t-1}(u) = \frac{-\sum_{j=1}^{p} \dot{X}_{t-j}^2(u)}{[1 + \sum_{j=1}^{p} X_{t-j}(u)^2]^2} \begin{pmatrix} 1 \\ X_{t-1}(u)^2 \\ \cdots \\ X_{t-p}(u)^2 \end{pmatrix} + \frac{1}{1 + \sum_{j=1}^{p} X_{t-j}(u)^2} \begin{pmatrix} 0 \\ \dot{X}_{t-1}^2(u) \\ \cdots \\ \dot{X}_{t-p}^2(u) \end{pmatrix}.$$

It is interesting to note that, like $\{\mathcal{M}_t(u)\}_t$, $\{\dot{\mathcal{M}}_t(u)\}_t$ is a sequence of vector martingale differences. We will use $\dot{\mathcal{M}}_t(u)$ to refine the approximation $\mathcal{R}_{t_0,N}(u_0)$ and show $\mathcal{R}_{t_0,N}(u_0) \approx \tilde{\mathcal{R}}_{t_0,N}(u_0)$, where

$$\tilde{\mathcal{R}}_{t_0,N}(u_0) = \sum_{k=0}^{t_0-1} \lambda(1 - \lambda F(u_0)) \left( \frac{t_0 - k}{N} - u_0 \right) \{\dot{\mathcal{M}}_{t_0-k}(u_0) + F(u_0)\underline{\dot{a}}(u_0)\}. \quad (65)$$

We use this to obtain Theorem 2.3.

**Lemma 5.1.** *Suppose Assumption 2.2 holds with some $r \geq 2$. Then we have*

$$|\dot{\mathcal{M}}_t(u) - \dot{\mathcal{M}}_t(v)|_1 \leq K|u - v|^{\beta'}|Z_t^2 + 1|L_t, \quad (66)$$

*where $L_t = \{1 + \sum_{k=1}^{p} W_{t-k}\}^2$, with $(\mathbb{E}(L_t)^{r/2})^{2/r} < \infty$, and almost surely*

$$\mathcal{M}_t(v) = \mathcal{M}_t(u) + (v - u)\dot{\mathcal{M}}_t(u) + (u - v)^{1+\beta'} R_t(u, v), \quad (67)$$

*where $|R_t(u, v)| \leq L_t$.*

**Proof.** By using (63), $\underline{a}(u) \in \text{Lip}(1 + \beta')$ and the fact that $|F_{t-1}(u)|_1 \leq (p+1)^2$ we have

$$|\dot{\mathcal{M}}_t(u) - \dot{\mathcal{M}}_t(v)|_1$$
$$\leq K|Z_t^2 + 1|\Big\{|F_{t-1}(u) - F_{t-1}(v)|_1$$
$$\qquad + |\dot{F}_{t-1}(u) - \dot{F}_{t-1}(v)|_1 + \sup_u |\dot{F}_{t-1}(u)|_1|u - v|\Big\}. \quad (68)$$

In order to bound (68), we consider $F_t(u)$ and its derivatives. We see that $|F_{t-1}(u) - F_{t-1}(v)|_1 \leq K|u - v|L_t$, thus bounding the first term on the right-hand side of (68). To obtain the other bounds we use (64). Now by using (22) and Lemma 4.2 we have that $\sup_u |\dot{F}_{t-1}(u)|_1 \leq KL_t$ and $|\dot{F}_{t-1}(u) - \dot{F}_{t-1}(v)|_1 \leq KL_t^2|u - v|^{\beta}$. Altogether this verifies (66). By using the Cauchy–Schwartz inequality we have that $(\mathbb{E}(L_t)^{r/2})^{2/r} < \infty$.

Finally, we prove (68). Since $L_t$ is a well-defined random variable almost surely all the paths of $\mathcal{M}_t(u) \in \text{Lip}(1 + \beta')$. Therefore, there exists a set $\mathcal{N}$ such that $P(\mathcal{N}) = 0$, and for every $\omega \in \mathcal{N}^c$ we can make a Taylor expansion of $\mathcal{M}_t(u, \omega)$ about $w$ and obtain (67). □

**Lemma 5.2.** *Suppose Assumption* 2.2 *holds with* $r \geq 2$. *Then if* $|t_0/N - u_0| < 1/N$ *we have*

$$\mathcal{R}_{t_0,N}(u_0) = \tilde{\mathcal{R}}_{t_0,N}(u_0) + R_N^{(4)}, \tag{69}$$

*where* $\mathbb{E}(|R_N^{(4)}|^r)^{1/r} = O(1/(N\lambda)^{1+\beta'})$.

**Proof.** To obtain the result we find the Taylor expansion of $\underline{a}((t-k)/N)$ and $\mathcal{M}_{t-k}((t-k)/N)$ about $u_0$, and substitute this into $\mathcal{R}_{t_0,N}$. Using Lemma 5.1, we obtain the desired result. $\qquad\square$

In the next lemma we consider the mean and variance of the bias and stochastic terms $\mathcal{R}_{t_0,N}(u_0)$ and $\mathcal{L}_{t_0}(u_0)$, which we use to obtain an asymptotic expression for the bias. We will use the following results. Since $\inf_u \lambda_{\min}(F(u)) > C > 0$ (see (97)), we have

$$\sum_{k=0}^{t-1} \lambda \{I - \lambda F(u)\}^k F(u) \to I, \tag{70}$$

$$\left\| \sum_{k=0}^{t-1} \lambda^2 \{I - \lambda F(u)\}^{2k} \left(\frac{k}{N}\right)^2 \right\| = O\left(\frac{1}{N^2\lambda}\right), \tag{71}$$

*if* $\lambda \to 0$, $t\lambda \to \infty$ *as* $t \to \infty$.

**Lemma 5.3.** *Suppose Assumption* 2.2 *holds with* $r > 4$. *Let* $\mathcal{R}_{t_0,N}(u_0)$, $\mathcal{L}_{t_0}(u_0)$ *and* $\Sigma(u)$ *be defined as in* (8) *and* (13), *respectively. Then if* $|t_0/N - u_0| < 1/N$ *we have*

$$\mathbb{E}(\mathcal{R}_{t_0,N}(u_0)) = -\frac{1}{N\lambda} F(u_0)^{-1} \underline{\dot{a}}(u) + O\left(\frac{1}{(N\lambda)^{1+\beta'}}\right), \tag{72}$$

$$\mathrm{var}(\mathcal{R}_{t_0,N}(u_0)) = O\left(\frac{1}{N^2\lambda} + \frac{1}{N}\right), \tag{73}$$

$$\mathbb{E}(\mathcal{L}_{t_0}(u_0)) = 0 \quad and \quad \mathrm{var}(\mathcal{L}_{t_0}(u_0)) = \lambda\Sigma(u_0) + o(\lambda). \tag{74}$$

**Proof.** Since $\{\partial \mathcal{M}_k(u)/\partial u\}$ are martingale differences, by applying (69) to (70) we have (72).

We now show (73). By using (71), $\sup_u |\partial \underline{a}(u)/\partial u| < \infty$ and $\sup_u |F_t(u)| < (p+1)^2$, we have

$$\mathrm{var}\{\mathcal{R}_{t_0,N}(u_0)\}$$
$$= \mu_4 \sum_{k=0}^{t_0-p-1} \lambda^2 \{I - \lambda F(u_0)\}^{2k} \left(\frac{k}{N}\right)^2$$
$$\qquad\qquad \times \mathrm{var}\{F_{t_0-k-1}(u_0)\underline{\dot{a}}(u_0) + \dot{F}_{t_0-k-1}(u_0)\underline{a}(u_0)\} + O\left(\frac{1}{N^2}\right)$$

$$= O\left(\frac{1}{N^2\lambda} + \frac{1}{N^2}\right).$$

It is straightforward to show (74) using (100).

Finally, since $\mathcal{L}_{t_0}(u_0)$ is the sum of martingale differences, $\mathbb{E}(\mathcal{L}_{t_0}(u_0)) = 0$.   □

From the above lemma it immediately follows that

$$\mathcal{R}_{t_0,N}(u_0) = -\frac{1}{N\lambda}F(u_0)^{-1}\dot{\underline{a}}(u_0) + O_p\left\{\frac{1}{(N\lambda)^{1+\beta'}} + \frac{1}{N\lambda^{1/2}}\right\} \tag{75}$$

and $(N\lambda)\mathcal{R}_{t_0,N}(u_0) \xrightarrow{\ell_2} -F(u_0)^{-1}\dot{\underline{a}}(u_0)$. We now use this prove Theorem 2.3.

**Proof of Theorem 2.3.** By substituting (75) into (60) (with $\beta = 1$) we have

$$\hat{\underline{a}}_{t_0,N} - \underline{a}(u_0) = \frac{-1}{N\lambda}F(u_0)^{-1}\dot{\underline{a}}(u_0) + \mathcal{L}_{t_0}(u_0) + \delta'_N R_N^{(5)}, \tag{76}$$

where $(\mathbb{E}|R_N^{(5)}|^2)^{1/2} < \infty$ and

$$\delta'_N = \frac{1}{(N\lambda)^{1+\beta'}} + \frac{1}{(N\lambda)^2} + \frac{1}{\lambda^{1/2}N} + \lambda + \frac{1}{N^{1/2}}.$$

By using the above and $\text{var}(\mathcal{L}_{t_0}(u_0)) = \lambda\Sigma(u_0)$ we have

$$\mathbb{E}|\hat{\underline{a}}_{t_0,N} - \underline{a}(u_0)|^2 = \lambda\,\text{tr}\{\Sigma(u_0)\} + \frac{1}{(N\lambda)^2}|F(u_0)^{-1}\dot{\underline{a}}(u_0)|^2 + O([(N\lambda)^{-1} + \sqrt{\lambda} + \delta'_N]\delta'_N),$$

which gives (15). In order to prove (16) we use (76) and Proposition 4.4. We first note that if $\lambda^{-1/2}/(N\lambda)^{1+\beta'} \to 0$, then by using (76) we have

$$\lambda^{-1/2}\{\hat{\underline{a}}_{t_0,N} - \underline{a}(u_0)\} + \lambda^{-1/2}\frac{1}{\lambda N}F(u_0)^{-1}\dot{\underline{a}}(u_0) = \lambda^{-1/2}\mathcal{L}_{t_0}(u_0) + o_p(1).$$

Therefore by using Proposition 4.4 we have (16).   □

# Appendix

## A.1. Mixingale properties of $\phi(\mathcal{X}_{t,N})$

Our object in this section is to prove Proposition 4.1. We do this by using the random vector representation of the tvARCH process $\{\mathcal{X}_{t,N}\}_t$. Let

$$A_t(u) = \begin{pmatrix} 0 & 0 & 0 & 0 & 0 & 0 & 0 \\ 0 & a_1(u)Z_t^2 & a_2(u)Z_t^2 & a_3(u)Z_t^2 & \dots & a_{p-1}(u)Z_t^2 & a_p(u)Z_t^2 \\ 0 & 1 & 0 & 0 & \dots & 0 & 0 \\ 0 & 0 & 1 & 0 & \dots & \dots & 0 \\ \vdots & \vdots & \vdots & \ddots & \ddots & \vdots & \vdots \\ \vdots & \vdots & \vdots & \ddots & \ddots & \vdots & \vdots \\ 0 & 0 & 0 & 0 & \dots & 1 & 0 \end{pmatrix}, \quad (77)$$

$\underline{b}_t(u)^\top = (1, a_0(u)Z_t^2, \underline{0}_{p-1}^\top)$. By using the definition of the tvARCH process given in (1) we have that the tvARCH vectors $\{\mathcal{X}_{t,N}\}_t$ satisfy the representation

$$\mathcal{X}_{t,N} = A_t\left(\frac{t}{N}\right)\mathcal{X}_{t-1,N} + \underline{b}_t\left(\frac{t}{N}\right). \tag{78}$$

Equation (78) looks like a vector autoregressive process; the difference is that $A_t(t/N)$ is a random matrix. However, similar to the vector autoregressive case, it can be shown that the product $\prod_{k=0}^{t} A_k(k/N)$ decays exponentially. It is this property which we use to prove Proposition 4.1.

Let $\mathcal{A}_N(t,j) = \{\prod_{i=0}^{j-1} A_{t-i}((t-i)/N)\}$, $\mathcal{A}(u,t,j) = \{\prod_{i=0}^{j-1} A_{t-i}(u)\}$, $\mathcal{A}_N(t,0) = I_{p+1}$ and $\mathcal{A}(u,t,0) = I_{p+1}$. By expanding (78) we have

$$\mathcal{X}_{t,N} = \mathcal{A}_N(t,k)\mathcal{X}_{t-k,N} + \sum_{j=0}^{k-1} \mathcal{A}_N(t,j)\underline{b}_{t-j}\left(\frac{t-j}{N}\right). \tag{79}$$

Suppose $A$ is a $n \times m$ dimensional random matrix, with $(i,j)$th element $A_{ij}$, and define the $n \times m$ dimensional matrix $[A]_q$ where $[A]_q = \{(\mathbb{E}|A_{ij}^q|)^{1/q}; i = 1,\dots,n, j = 1,\dots,m\}$. Now by using Proposition 2.1 in Subba Rao [14] and Corollary A.2 below, we have

$$\|[\mathcal{A}_N(t,k)]_q\| \le K\rho^k \quad \text{and} \quad \|[\mathcal{A}(u,t,k)]_q\| \le K\rho^k. \tag{80}$$

**Proof of Proposition 4.1.** We first prove (24). Let $\mathcal{C}_N(t,k) := \sum_{j=0}^{k-1} \mathcal{A}_N(t,j)\underline{b}_{t-j}((t-j)/N)$, that is,

$$\mathcal{X}_{t,N} = \mathcal{A}_N(t,k)\mathcal{X}_{t-k,N} + \mathcal{C}_N(t,k), \tag{81}$$

with $\mathcal{A}_N(t,k), \mathcal{C}_N(t,k) \in \sigma(Z_t,\dots,Z_{t-k+1})$ and $\mathcal{X}_{t-k,N} \in \mathcal{F}_{t-k}$. In particular, we have

$$\mathbb{E}\{\phi(\mathcal{C}_N(t,k))|\mathcal{F}_{t-k}\} = \mathbb{E}\{\phi(\mathcal{C}_N(t,k))\}. \tag{82}$$

Furthermore, by using Minkowski's inequality it can be shown that

$$\{\mathbb{E}(|\mathcal{A}_N(t,k)\mathcal{X}_{t-k,N}|^q|\mathcal{F}_{t-k})\}^{1/q} \leq K|[\mathcal{A}_N(t,k)]_q\mathcal{X}_{t-k,N}|_1.$$

The Lipschitz continuity of $\phi$ and (79) now imply $|\phi(\mathcal{X}_{t,N}) - \phi(\mathcal{C}_N(t,k))| \leq K|\mathcal{A}_N(t,k)\mathcal{X}_{t-k,N}|$. Therefore, by using the above we obtain

$$\begin{aligned}
&|\mathbb{E}\{\phi(\mathcal{X}_{t,N})|\mathcal{F}_{t-k})\} - \mathbb{E}\{\phi(\mathcal{X}_{t,N})\}|_1 \\
&= |\mathbb{E}_{t-k}\{\phi(\mathcal{X}_{t,N}) - \phi(\mathcal{C}_N(t,k))|\mathcal{F}_{t-k})\} - \mathbb{E}\{\phi(\mathcal{X}_{t,N}) - \phi(\mathcal{C}_N(t,k))\}|_1 \\
&\leq K(|[\mathcal{A}_N(t,k)]_1\mathcal{X}_{t-k,N}|_1 + \mathbb{E}\{|[\mathcal{A}_N(t,k)]_1\mathcal{X}_{t-k,N}|_1\}) \\
&\leq K\{\||[\mathcal{A}_N(t,k)]_1\||\mathcal{X}_{t-k,N}|_1 + \mathbb{E}(\||[\mathcal{A}_N(t,k)]_1\||\mathcal{X}_{t-k,N}|_1)\} \\
&\leq K\rho^k(1 + |\mathcal{X}_{t-k,N}|_1)
\end{aligned} \tag{83}$$

since $\mathbb{E}|\mathcal{X}_{t-k,N}|_1$ is uniformly bounded, thus giving (24). The proof of (25) is the same as the proof above, so we omit details. Inequality (26) follows from (24) with the triangular inequality.

To prove (27) we use (81). Since

$$\mathbb{E}\{\phi(\mathcal{C}_N(t,k))|\mathcal{F}_{t-k}\} = \mathbb{E}\{\phi(\mathcal{C}_N(t,k))\} = \mathbb{E}\{\phi(\mathcal{C}_N(t,k))|\mathcal{F}_{t-k-j}\},$$

we obtain, as in (83),

$$\begin{aligned}
&|\mathbb{E}\{\phi(\mathcal{X}_{t,N})|\mathcal{F}_{t-k}\} - \mathbb{E}\{\phi(\mathcal{X}_{t,N})|\mathcal{F}_{t-k-j}\}|_1 \\
&\leq K\mathbb{E}\{|[\mathcal{A}_N(t,k)]_1\mathcal{X}_{t-k,N}|_1|\mathcal{F}_{t-k}\} + \mathbb{E}\{|[\mathcal{A}_N(t,k)]_1\mathcal{X}_{t-k,N}|_1|\mathcal{F}_{t-k-j}\} \\
&\leq K\rho^k(|\mathcal{X}_{t-k,N}|_1 + |\mathcal{X}_{t-k,N}|_1),
\end{aligned}$$

which gives (27).

To prove (28) we will use (83). We first note that by using Minkowski's inequality and the equivalence of norms, there exists a constant $K$ independent of $\mathcal{X}_{t,N}$ such that

$$\{\mathbb{E}(|\mathcal{X}_{t,N}|^q|\mathcal{F}_{t-k})\}^{1/q} \leq \{\mathbb{E}(|\mathcal{A}_N(t,k)\mathcal{X}_{t-k,N}|^q|\mathcal{F}_{t-k})\}^{1/q} + \{\mathbb{E}(|C_N(t,k)|^q|\mathcal{F}_{t-k})\}^{1/q}.$$

Now by using $\{\mathbb{E}(|\mathcal{A}_N(t,k)\mathcal{X}_{t-k,N}|^q|\mathcal{F}_{t-k})\}^{1/q} \leq K\rho^k|\mathcal{X}_{t-k,N}|_1$ and

$$[\mathbb{E}|C_N(t,k)|^q|\mathcal{F}_{t-k}]^{1/q} \leq K\sum_{j=0}^{k-1}\left|[\mathcal{A}_N(t,j)]_q\left[\underline{b}_{t-j}\left(\frac{t-j}{N}\right)\right]_q\right|_1 \leq \frac{K}{1-\rho},$$

we have

$$\mathbb{E}\{|\mathcal{X}_{t,N}|^q|\mathcal{F}_{t-k}\} \leq \left\{K\rho^k|\mathcal{X}_{t-k,N}|_1 + \frac{K}{1-\rho}\right\}^q \leq K|\mathcal{X}_{t-k,N}|_1^q,$$

hence we obtain (28). □

We use the corollary below to prove Lemma A.7 at the end of Section A.3.

**Corollary A.1.** *Suppose Assumption* 2.1 *holds with* $r = q_0$ *and* $\mathbb{E}(Z_0^{4q_0}) < \infty$. *Let* $f, g : \mathbb{R}^{p+1} \to \mathbb{R}^{(p+1) \times (p+1)}$ *be Lipschitz continuous functions, such that for all positive* $\underline{x} \in \mathbb{R}^{p+1}$, $|f(\underline{x})|_{\mathrm{abs}} \leq \mathbb{I}_{p+1}$ *and* $|g(\underline{x})|_{\mathrm{abs}} \leq \mathbb{I}_{p+1}$, *where* $\mathbb{I}_{p+1}$ *is a* $(p+1) \times (p+1)$ *matrix with one in all the elements, and* $|A|_{\mathrm{abs}} = (|A_{i,j}|: i, j = 1, \ldots, p+1)$. *Then for* $q \leq q_0$ *we have*

$$
\begin{aligned}
&[\mathbb{E}(\mathbb{E}\{(Z_{t-k+1}^2 - 1)f(\mathcal{X}_{t-k,N})g(\mathcal{X}_t(u))|\mathcal{F}_{t-k-j}\} \\
&\qquad - \mathbb{E}\{(Z_{t-k+1}^2 - 1)f(\mathcal{X}_{t-k,N})g(\mathcal{X}_t(u))\})^q]^{1/q} \\
&\leq K\{\mathbb{E}(Z_0^4) + 1\}\rho^j((\mathbb{E}|\mathcal{X}_{t-k-j,N}|^q)^{1/q} + (\mathbb{E}|\mathcal{X}_{t-k-j}(u)|^q)^{1/q})
\end{aligned}
\tag{84}
$$

*and*

$$
\begin{aligned}
&[\mathbb{E}(\mathbb{E}\{(Z_{t-k+1}^2 - 1)f(\mathcal{X}_{t-k}(u))g(\mathcal{X}_t(u))|\mathcal{F}_{t-k-j}\} \\
&\qquad - \mathbb{E}\{(Z_{t-k+1}^2 - 1)f(\mathcal{X}_{t-k}(u))g(\mathcal{X}_t(u))\})^q]^{1/q} \\
&\leq K\{\mathbb{E}(Z_0^4) + 1\}\rho^j((\mathbb{E}|\mathcal{X}_{t-k-j}(u)|^q)^{1/q} + (\mathbb{E}|\mathcal{X}_{t-k-j}(u)|^q)^{1/q}),
\end{aligned}
\tag{85}
$$

*for* $j, k \geq 0$, *where* $\rho$ *is such that* $1 - \eta < \rho < 1$.

**Proof.** We give the proof of (84) only; that of (85) is the same. We use the notation introduced in the proof of Proposition 4.1 and let $\mathcal{C}(u, t, k) = \sum_{j=0}^{k-1} \mathcal{A}(u, t, j)\underline{b}_{t-j}(u)$, that is, $\mathcal{X}_t(u) = \mathcal{A}(u, t, k+j)\mathcal{X}_{t-k-j}(u) + \mathcal{C}(u, t, k+j)$.

We use (82). Since $|f(\underline{x})|_{\mathrm{abs}} \leq \mathbb{I}_{p+1}$ and $|g(\underline{x})|_{\mathrm{abs}} \leq \mathbb{I}_{p+1}$ we have

$$
\begin{aligned}
&|(Z_{t-k+1}^2 - 1)(f(\mathcal{X}_{t-k,N})g(\mathcal{X}_t(u)) - f(\mathcal{C}_N(t-k,j))g(\mathcal{C}(u, t, k+j)))| \\
&\leq |Z_{t-k+1}^2 + 1|(|\mathcal{A}_N(t-k, j)\mathcal{X}_{t-k-j,N}|_1 + |\mathcal{A}(u, t, k+j)\mathcal{X}_{t-k-j}(u)|_1).
\end{aligned}
\tag{86}
$$

Since $\mathcal{A}(u, t, k-1)$, $(Z_{t-k}^2 + 1)A_{t-k-1}(u)$ and $\mathcal{A}(u, t-k, j)$ are independent matrices and by using (80), we have

$$
\begin{aligned}
&\mathbb{E}\|[(Z_{t-k+1}^2 + 1)\mathcal{A}(u, t, k+j)]_1\| \\
&\leq \mathbb{E}\|[\mathcal{A}(u, t, k-1)]_1\|\|\|[(Z_{t-k+1}^2 + 1)A_{t-k+1}(u)]_1\|_1\|[\mathcal{A}(u, t-k, j)]_1\| \\
&\leq K\mathbb{E}(Z_{t-k+1}^4 + 1)\rho^{k+j-1}.
\end{aligned}
\tag{87}
$$

Considering the conditional expectation of (86), by using (87) and $|\mathbb{E}_{t-k-j}\{(Z_{t-k+1}^2 - 1)\mathcal{A}_N(t-k, j)\mathcal{X}_{t-k-j,N}\}| \leq K\rho^j|\mathcal{X}_{t-k-j,N}|_1$, we have

$$
\begin{aligned}
&|\mathbb{E}[(Z_{t-k+1}^2 - 1)\{f(\mathcal{X}_{t-k,N})g(\mathcal{X}_t(u)) - f(\mathcal{C}_N(t-k,j))g(\mathcal{C}(u, t, k+j))\}|\mathcal{F}_{t-k-1}]| \\
&\leq K\rho^j|\mathcal{X}_{t-k-j,N}|_1\mathbb{E}(Z_{t-k+1}^2 + 1) + K\mathbb{E}(Z_{t-k+1}^4 + 1)\rho^{k+j-1}|\mathcal{X}_{t-k-j}(u)|_1 \\
&\leq K\mathbb{E}(Z_0^4 + 1)\rho^j(|\mathcal{X}_{t-k-j,N}|_1 + |\mathcal{X}_{t-k-j}(u)|_1).
\end{aligned}
$$

Using similar arguments, we obtain the bound

$$(\mathbb{E}|\mathbb{E}[(Z_{t-k+1}^2 - 1)\{f(\mathcal{X}_{t-k,N})g(\mathcal{X}_t(u)) - f(\mathcal{C}_N(t-k,j))g(\mathcal{C}(u,t,k+j))\}]|^r)^{1/r}$$

$$\leq K\mathbb{E}(Z_0^4 + 1)\rho^j (\mathbb{E}(|\mathcal{X}_{t-k-j,N}|^r)^{1/r} + (\mathbb{E}|\mathcal{X}_{t-k-j}(u)|^r)^{1/r}),$$

leading to the result. $\qquad\qquad\square$

## A.2. Persistence of excitation

As we mentioned in Section 1, a crucial component in the analysis of many stochastic, recursive algorithms is to establish persistence of excitation of the transition random matrix in the algorithm: in our case this implies showing Theorem 4.1. Intuitively, it is clear that the verification of this result depends on showing that the smallest eigenvalue of the conditional expectation of of the semi-positive definite matrix $\mathcal{X}_{t,N}\mathcal{X}_{t,N}^\top/\|\mathcal{X}_{t,N}\|_1^2$ is bounded away from zero. In particular this is one of the major conditions given in Moulines *et al.* ([11], Theorem 16), where conditions were given under which persistence of excitation can be established. In fact in this section we verify the conditions in Moulines *et al.* ([11], Theorem 16) to prove Theorem 4.1.

Suppose $X$ is a random variable and define $\mathbb{E}_t(X) = \mathbb{E}(X|\mathcal{F}_t)$.

**Lemma A.1.** *Suppose that Assumption 2.1 holds with $r = 2$. Then we have*

$$\lambda_{\min}(\mathbb{E}\{\mathcal{X}_t(u)\mathcal{X}_t(u)^\top|\mathcal{F}_{t-k}\}) > C \tag{88}$$

*and, for $N$ large enough,*

$$\lambda_{\min}(\mathbb{E}\{\mathcal{X}_{t,N}\mathcal{X}_{t,N}^\top|\mathcal{F}_{t-k}\}) > \frac{C}{2}, \tag{89}$$

*for $k \geq (p+1)$, where $C$ is a finite constant independent of $t, N$ and $u$.*

**Proof.** We will prove (89); the proof of (88) is similar. We partition $\mathcal{X}_{t,N}\mathcal{X}_{t,N}^\top$ as

$$\mathcal{X}_{t,N}\mathcal{X}_{t,N}^\top = \Delta_1 + \Delta_2, \tag{90}$$

where $\Delta_1$ and $\Delta_2$ are positive definite matrices, $\Delta_1 \in \sigma(Z_t, \dots, Z_{t-p})$ and $\Delta_2 \in \mathcal{F}_t$. This implies $\lambda_{\min}\{\mathbb{E}_{t-k}(\mathcal{X}_{t,N}\mathcal{X}_{t,N}^\top)\} \geq \lambda_{\min}\{\mathbb{E}(\Delta_1)\}$, if $k \geq p+1$, which allows us to obtain a uniform lower bound for $\lambda_{\min}\{\mathbb{E}_{t-k}(\mathcal{X}_{t,N}\mathcal{X}_{t,N}^\top)\}$. To facilitate this we represent $\mathcal{X}_{t,N}$ in terms of martingale vectors. By using (1) we have

$$X_{t,N}^2 = a_0\left(\frac{t}{N}\right) + \sum_{i=1}^p a_i\left(\frac{t}{N}\right)X_{t-i,N}^2 + (Z_t^2 - 1)\sigma_{t,N}^2$$

and

$$\mathcal{X}_{t,N} = \Theta\left(\frac{t}{N}\right)\mathcal{X}_{t-1,N} + (Z_t^2 - 1)\sigma_{t,N}^2 D, \tag{91}$$

where $D$ is a $(p+1)$-dimensional vector with $D_i = 0$ if $i \neq 2$ and $D_2 = 1$, and

$$\Theta(u) = \begin{pmatrix} 1 & 0 & 0 & 0 & 0 & 0 \\ a_0(u) & a_1(u) & a_2(u) & \dots & a_{p-1}(u) & a_p(u) \\ 0 & 1 & 0 & \dots & 0 & 0 \\ 0 & 0 & 1 & \dots & 0 & 0 \\ \vdots & \vdots & \ddots & \vdots & \vdots & \\ 0 & 0 & 0 & \dots & 1 & 0 \end{pmatrix}.$$

With $p+1$ iterations of (91)

$$\mathcal{X}_{t,N} = \sum_{i=0}^{p}(Z_{t-i}^2 - 1)\sigma_{t-i,N}^2 \left\{\prod_{j=0}^{i-1}\Theta\left(\frac{t-j}{N}\right)\right\}D + \left\{\prod_{j=0}^{p}\Theta\left(\frac{t-j}{N}\right)\right\}\mathcal{X}_{t-p-1,N}.$$

Let $\mu_4 = \mathbb{E}(Z_t^2 - 1)^2$, $\Theta_{t,N}(i) = \prod_{j=0}^{i-1}\Theta((t-j)/N)$ and $B_{t,N}(i) = \Theta_{t,N}(i)DD^\top\Theta_{t,N}(i)^\top$. Since $\{(Z_t^2 - 1)\sigma_{t,N}^2\}_t$ are martingale differences, we have

$$\mathbb{E}(\mathcal{X}_{t,N}\mathcal{X}_{t,N}^\top | \mathcal{F}_{t-k})$$

$$= \mu_4 \sum_{i=0}^{p} \mathbb{E}_{t-k}(\sigma_{t-i,N}^4)\Theta_{t,N}(i)DD^\top\Theta_{t,N}(i)^\top$$

$$+ \Theta_{t,N}(p+1)\mathbb{E}_{t-k}(\mathcal{X}_{t-p-1,N}\mathcal{X}_{t-p-1,N}^\top)\Theta_{t,N}(p+1)^\top \qquad \text{for } k \geq p+1. \tag{92}$$

Since the matrices above are non-negative definite, we have that

$$\lambda_{\min}\{\mathbb{E}(\mathcal{X}_{t,N}\mathcal{X}_{t,N}^\top | \mathcal{F}_{t-k})\} \geq \lambda_{\min}\left\{\mu_4 \sum_{i=0}^{p}\mathbb{E}_{t-k}(\sigma_{t-i,N}^4)B_{t,N}(i)\right\}. \tag{93}$$

We now refine $\mu_4 \sum_{i=0}^{p}\mathbb{E}_{t-k}(\sigma_{t-i,N}^4)B_{t,N}(i)$ to obtain $\Delta_1$. By using (79) we have

$$\mathcal{X}_{t,N} = \sum_{j=0}^{p}\mathcal{A}_N(t,j)\underline{b}_{t-j}\left(\frac{t-j}{N}\right) + \mathcal{A}_N(t,p+1)\mathcal{X}_{t-p-1,N}.$$

Therefore, by using (1) and $X_{t-i,N}^2 = (\mathcal{X}_{t,N})_{i+1}$, we have

$$\sigma_{t-i,N}^2 = \frac{1}{Z_{t-i}^2}\sum_{j=0}^{p}\left\{\mathcal{A}_N(t,j)\underline{b}_{t-j}\left(\frac{t-j}{N}\right)\right\}_{i+1} + \frac{1}{Z_{t-i}^2}\{\mathcal{A}_N(t,p+1)\mathcal{X}_{t-p-1,N}\}_{i+1}$$

$$= H_{t-i,N}(t) + G_{t-i,N}(t), \tag{94}$$

where $H_{t-i,N}(t) \in \sigma(Z_{t-i}, \ldots, Z_{t-p})$ and $G_{t-i,N}(t) \in \mathcal{F}_{t-i}$. Since $H_{t-i,N}(t)$ and $G_{t-i,N}(t)$ are positive this implies, with (93),

$$\lambda_{\min}\{\mathbb{E}(\mathcal{X}_{t,N}\mathcal{X}_{t,N}^\top|\mathcal{F}_{t-k})\} \geq \lambda_{\min}\{P_{t,N}\},$$

where

$$P_{t,N} = \mu_4 \sum_{i=0}^{p} \mathbb{E}(H_{t-i,N}(t)^2) B_{t,N}(i).$$

To bound this we define the corresponding terms for the stationary approximation $X_t(u)$. We set

$$P(u) = \mu_4 \sum_{i=0}^{p} \mathbb{E}(H_{t-i}(u,t)^2) B(u,i),$$

where

$$H_{t-i}(u,t) = \frac{1}{Z_{t-i}^2} \sum_{j=0}^{p} \{\mathcal{A}(u,t,j)\underline{b}_{t-j}(u)\}_{i+1}$$

and $B(u,i) = \Theta(u)^i DD^\top(\Theta(u)^i)^\top$. A close inspection of the above calculation steps reveals that

$$P(u) = \mathbb{E}\{\underline{Y}_p(u)\underline{Y}_p(u)^\top\}$$

with $\underline{Y}_p(u)$ from Assumption 2.1(iv). Therefore,

$$\lambda_{\min}\left\{P\left(\frac{t}{N}\right)\right\} \geq \inf_u \lambda_{\min}\{P(u)\} \geq C.$$

Since $\{a_j(\cdot)\}_j$ is $\beta$-Lipschitz, we have, for $i = 0, \ldots, p-1$,

$$\left|\mathbb{E}\left(H_{t-i}\left(\frac{t}{N},t\right)^2\right) - \mathbb{E}(H_{t-i,N}(t)^2)\right|_1 \leq \frac{K}{N^\beta},$$

$$\left|B_{t,N}(i) - B\left(\frac{t}{N},i\right)\right|_1 \leq \frac{K}{N^\beta}.$$

Therefore $\|P_{t,N} - P(t/N)\| \leq |P_{t,N} - P(t/N)|_1 \leq K/N^\beta$, which leads to

$$\lambda_{\min}\{P_{t,N}\} = \lambda_{\min}\left\{P\left(\frac{t}{N}\right) + \left[P_{t,N} - P\left(\frac{t}{N}\right)\right]\right\}$$

$$\geq \lambda_{\min}\left\{P\left(\frac{t}{N}\right)\right\} - \left\|P_{t,N} - P\left(\frac{t}{N}\right)\right\| \geq C - \frac{K}{N^\beta}$$

and therefore to

$$\lambda_{\min}(\mathbb{E}\{\mathcal{X}_{t,N}\mathcal{X}_{t,N}^{\top}|\mathcal{F}_{t-k,N}\}) > C - \frac{K}{N^{\beta}}. \tag{95}$$

Thus for $N$ large enough we have (89). $\qquad\square$

**Lemma A.2.** *Suppose that Assumption* 2.1 *holds with* $r = 4$. *Then, for* $k \geq p+1$, *we have for* $N$ *large enough,*

$$\lambda_{\min}\left(\mathbb{E}\left\{\frac{\mathcal{X}_{t,N}\mathcal{X}_{t,N}^{\top}}{|\mathcal{X}_{t,N}|_1^2}\Big|\mathcal{F}_{t-k}\right\}\right) > \frac{C}{|\mathcal{X}_{t-k,N}|_1^4} \tag{96}$$

*and*

$$\lambda_{\min}\left(\mathbb{E}\left\{\frac{\mathcal{X}_t(u)\mathcal{X}_t(u)^{\top}}{|\mathcal{X}_t(u)|_1^2}\right\}\right) > C, \tag{97}$$

*where* $C$ *is a constant independent of* $t, N$ *and* $u$.

**Proof.** We first prove (96). By definition,

$$\lambda_{\min}\left\{\mathbb{E}_{t-k}\left(\frac{\mathcal{X}_{t,N}\mathcal{X}_{t,N}^{\top}}{|\mathcal{X}_{t,N}|_1^2}\right)\right\} = \inf_{|\underline{x}|=1}\mathbb{E}_{t-k}\left(\frac{\underline{x}^{\top}\mathcal{X}_{t,N}}{|\mathcal{X}_{t,N}|_1}\right)^2.$$

Since

$$(\underline{x}^{\top}\mathcal{X}_{t,N})^2 = \frac{\underline{x}^{\top}\mathcal{X}_{t,N}}{|\mathcal{X}_{t,N}|_1}\underline{x}^{\top}\mathcal{X}_{t,N}|\mathcal{X}_{t,N}|_1$$

and $\sup_{|\underline{x}|=1}|\underline{x}^{\top}\mathcal{X}_{t,N}|_1 \leq |\mathcal{X}_{t,N}|_1$, we obtain by using Cauchy's inequality and (28),

$$\mathbb{E}_{t-k}(\underline{x}^{\top}\mathcal{X}_{t,N})^2 \leq \left\{\mathbb{E}_{t-k}\left(\frac{\underline{x}^{\top}\mathcal{X}_{t,N}}{|\mathcal{X}_{t,N}|_1}\right)^2\right\}^{1/2}\{\mathbb{E}_{t-k}(\underline{x}^{\top}\mathcal{X}_{t,N}|\mathcal{X}_{t,N}|_1)^2\}^{1/2}$$

$$\leq \left\{\mathbb{E}_{t-k}\left(\frac{\underline{x}^{\top}\mathcal{X}_{t,N}}{|\mathcal{X}_{t,N}|_1}\right)^2\right\}^{1/2}\{K|\mathcal{X}_{t-k,N}|_1^4\}^{1/2}.$$

Therefore by using the above and Lemma A.1 for large $N$, we obtain

$$\inf_{|\underline{x}|=1}\mathbb{E}_{t-k}\left(\frac{\underline{x}^{\top}\mathcal{X}_{t,N}}{|\mathcal{X}_{t,N}|_1}\right)^2 \geq \inf_{|\underline{x}|=1}\frac{[\mathbb{E}_{t-k}(\underline{x}^{\top}\mathcal{X}_{t,N})^2]^2}{\mathbb{E}_{t-k}(|\mathcal{X}_{t,N}|_1^4)} \geq \frac{C}{|\mathcal{X}_{t-k,N}|_1^4},$$

where $C$ is a positive constant, thus giving (96).

To prove (97), we use (88) to obtain

$$\mathbb{E}(\mathcal{X}_t(u)\mathcal{X}_t(u)^{\top}) = \mathbb{E}_{-\infty}(\mathcal{X}_t(u)\mathcal{X}_t(u)^{\top}) > C,$$

and using the arguments above we have

$$\lambda_{\min}\left(\mathbb{E}\left\{\frac{\mathcal{X}_t(u)\mathcal{X}_t(u)^\top}{|\mathcal{X}_t(u)|_1^2}\right\}\right) > \frac{C}{\mathbb{E}(|\mathcal{X}_{t-k}(u)|_1^4)}.$$

By Lemma 4.1 $\sup_u \mathbb{E}(|\mathcal{X}_t(u)|_1^4) < \infty$, which leads to (97). $\qquad\square$

**Corollary A.2.** *Suppose Assumption 2.1 holds with $r = 4$. Let $F(u)$ be defined as in (9). Then there exist $C$ and $\lambda_1$ such that, for all $\lambda \in [0, \lambda_1]$ and $u \in (0, 1]$, we have*

$$\lambda_{\max}\{I - \lambda F(u)\} \le 1 - \lambda C. \tag{98}$$

*There exist a $0 < \delta \le C$ and $K$ such that, for all $k$,*

$$\|\{I - \lambda F(u)\}^k\| \le K(1 - \lambda\delta)^k. \tag{99}$$

*Furthermore,*

$$\sum_{k=0}^{t-1} \lambda\{I - \lambda F(u)\}^{2k} \to \tfrac{1}{2}F(u)^{-1}, \tag{100}$$

*where $\lambda \to 0$ and $\lambda t \to \infty$ as $t \to \infty$.*

**Proof.** Inequality (98) follows directly from (97). Furthermore, since $(I - \lambda F(u))$ is symmetric matrix, we have $\|\{I - \lambda F(u)\}^k\| \le \|(I - \lambda F(u))\|^k \le (1 - \lambda\delta)^k$, that is, (99).

We now prove (100). We have

$$\lambda \sum_{k=0}^{t-1} \{I - \lambda F(u)\}^{2k} = \lambda\{I - (I - \lambda F(u))^2\}^{-1}\{I - (I - \lambda F(u))^{2t}\}.$$

Since $\lambda_{\min}\{F(u)\} > C$, for some $C > 0$, we have $|\{I - \lambda F(u)\}^{2t}|_1 \le \|\{I - \lambda F(u)\}^{2t}\|\|I_{p+1}| \to 0$ as $\lambda \to 0$, $\lambda t \to \infty$ and $t \to \infty$. Furthermore, $\lambda\{I - (I - \lambda F(u))^2\}^{-1} \to \tfrac{1}{2}F(u)^{-1}$. Together these give (100). $\qquad\square$

**Lemma A.3.** *Suppose that Assumption 2.1 holds with $r = 4$. For a sufficiently large $N$ and for every $R > 0$, there exist an $s_0 > p + 1$ and $C_1$ such that, for all $s \ge s_0$ and $t = 1, \ldots, N - s$, we have*

$$\mathbb{E}\left\{\lambda_{\min}\left(\sum_{k=t}^{t+s-1} \frac{\mathcal{X}_{k,N}\mathcal{X}_{k,N}^\top}{|\mathcal{X}_{k,N}|_1^2}\Big|\mathcal{F}_t\right)\right\} \ge C_1 I(|\mathcal{X}_{t,N}|_1 \le R),$$

*where $I$ denotes the identity function.*

**Proof.** The result can be proved using the methods given in Moulines *et al.* ([11], Lemma 19). We outline the proof. By using $\lambda_{\min}(A) \geq \lambda_{\min}(B) - \|A - B\|$ (see Moulines *et al.* [11], Lemma 19), we have

$$\mathbb{E}_t\left\{\lambda_{\min}\left(\sum_{k=t}^{t+s-1} \frac{\mathcal{X}_{k,N}\mathcal{X}_{k,N}^\top}{|\mathcal{X}_{k,N}|_1^2}\right)\right\}$$

$$\geq \sum_{k=t}^{t+s-1} \lambda_{\min}\left\{\mathbb{E}_t\left(\frac{\mathcal{X}_{k,N}\mathcal{X}_{k,N}^\top}{|\mathcal{X}_{k,N}|_1^2}\right)\right\} - \mathbb{E}_t\left\|\sum_{k=t}^{t+s-1}\left\{\frac{\mathcal{X}_{k,N}\mathcal{X}_{k,N}^\top}{|\mathcal{X}_{k,N}|_1^2} - \mathbb{E}_t\left(\frac{\mathcal{X}_{k,N}\mathcal{X}_{k,N}^\top}{|\mathcal{X}_{k,N}|_1^2}\right)\right\}\right\|.$$

We now evaluate an upper bound for the second term above. Let

$$\Delta_{k,N} = \frac{\mathcal{X}_{k,N}\mathcal{X}_{k,N}^\top}{|\mathcal{X}_{k,N}|_1^2} - \mathbb{E}_t\left(\frac{\mathcal{X}_{k,N}\mathcal{X}_{k,N}^\top}{|\mathcal{X}_{k,N}|_1^2}\right).$$

Then

$$\mathbb{E}_t\left\|\sum_{k=t}^{t+s-1}\Delta_{k,N}\right\|^2 \leq \mathbb{E}_t\left|\sum_{k=t}^{t+s-1}\Delta_{k,N}\right|^2 \leq \sum_{i,j=1}^{p+1}\sum_{k_1,k_2=t}^{t+s-1}\mathbb{E}_t((\Delta_{k_1,N})_{i,j}(\Delta_{k_2,N})_{i,j}), \quad (101)$$

and we require a bound on $\mathbb{E}_t((\Delta_{k_1,N})_{i,j}(\Delta_{k_2,N})_{i,j})$. For $k_1 = k_2 = k$ we obtain with (28),

$$\mathbb{E}_t\{(\Delta_{k,N})_{i,j}\}^2 \leq 2\mathbb{E}_t|\mathcal{X}_{k,N}|_1^4 \leq K|\mathcal{X}_{t,N}|_1^4. \quad (102)$$

Now let $k_1 \neq k_2$. Let $\phi(\underline{x}) = \underline{x}^\top\underline{x}/|\underline{x}|_1^2$, where $\underline{x} = (1, x_1, \ldots, x_p)$. Since $\phi \in \text{Lip}(1)$, by using (27) we have, for $k_1 < k_2$,

$$|\mathbb{E}_{k_1}((\Delta_{k_2,N})_{i,j})| = \left|\mathbb{E}_{k_1}\left(\frac{(\mathcal{X}_{k_2,N}\mathcal{X}_{k_2,N}^\top)_{i,j}}{|\mathcal{X}_{k_2,N}|_1^2}\right) - \mathbb{E}_t\left(\frac{(\mathcal{X}_{k_2,N}\mathcal{X}_{k_2,N}^\top)_{i,j}}{|\mathcal{X}_{k_2,N}|_1^2}\right)\right|$$

$$\leq K\rho^{k_2-k_1}(\mathbb{E}_t(|\mathcal{X}_{k_1,N}|_1) + |\mathcal{X}_{k_1,N}|_1). \quad (103)$$

By using (28), (102) and (103) we have

$$|\mathbb{E}_t\{(\Delta_{k_1,N})_{i,j}(\Delta_{k_2,N})_{i,j}\}| \leq |\mathbb{E}_t\{\mathbb{E}_{k_1}(\Delta_{k_1,N})_{i,j}(\Delta_{k_2,N})_{i,j}\}|$$

$$\leq |\mathbb{E}_t\{(\Delta_{k_1,N})_{i,j}\mathbb{E}_{k_1}(\Delta_{k_2,N})_{i,j}\}|$$

$$\leq \{\mathbb{E}_t(\Delta_{k_1,N})_{i,j}^2\}^{1/2}\{\mathbb{E}_t(\mathbb{E}_{k_1}(\Delta_{k_2,N})_{i,j})^2\}^{1/2}$$

$$\leq K\rho^{k_2-k_1}|\mathcal{X}_{t,N}|_1^3 \leq K\rho^{k_2-k_1}|\mathcal{X}_{t,N}|_1^4. \quad (104)$$

Substituting (104) into (101), we obtain $\mathbb{E}_t\|\sum_{k=t}^{t+s-1}\Delta_{k,N}\|^2 = Ks|\mathcal{X}_{t,N}|_1^4$, where $K$ is a constant independent of $s$. Therefore, by using this and (96), we have, for $N$ sufficiently large

$$\mathbb{E}_t\left\{\lambda_{\min}\left(\sum_{k=t}^{t+s-1} \frac{\mathcal{X}_{k,N}\mathcal{X}_{k,N}^\top}{|\mathcal{X}_{k,N}|_1^2}\right)\right\} \geq \frac{Cs}{2}|\mathcal{X}_{t,N}|_1^{-4} - K\{s|\mathcal{X}_{t,N}|_1^4\}^{1/2}.$$

Therefore for any $R > 0$,

$$\mathbb{E}_t \left\{ \lambda_{\min} \left( \sum_{k=t}^{t+s-1} \frac{\mathcal{X}_{k,N} \mathcal{X}_{k,N}^\top}{|\mathcal{X}_{k,N}|_1^2} \right) \right\} \geq \frac{s^{1/2}}{R^4} (C s^{1/2} - KR^2) I(|\mathcal{X}_{t,N}|_1 \leq R).$$

Now choose $s_0$ and a corresponding $C_1$ such that

$$C_1 = \frac{s_0^{1/2}}{R^4} (C s_0^{1/2} - KR^2) > 0.$$

Then it is clear that if $s > s_0$ then we have

$$\mathbb{E}_t \left\{ \lambda_{\min} \left( \sum_{k=t}^{t+s-1} \frac{\mathcal{X}_{k,N} \mathcal{X}_{k,N}^\top}{|\mathcal{X}_{k,N}|_1^2} \right) \right\} \geq \frac{s^{1/2}}{R^4} (C s^{1/2} - KR^2) I(|\mathcal{X}_{t,N}|_1 \leq R)$$

$$\geq C_1 I(|\mathcal{X}_{t,N}|_1 \leq R),$$

thus giving the result. $\square$

**Proof of Theorem 4.1.** We prove the result by verifying the conditions in Moulines *et al.* ([11], Theorem 16); see also Dahlhaus and Subba Rao [3]. Let $\phi_l := |\mathcal{X}_{l,N}|_1$, $V_l := |\mathcal{X}_{l,N}|_1$ and $A_l := F_{l,N} = \mathcal{X}_{l,N} \mathcal{X}_{l,N}^\top / |\mathcal{X}_{l,N}|_1^2$. We now show that conditions (a)–(d) in Moulines *et al.* ([11], Theorem 16) are satisfied. Let $k \geq p+1$ and $1 \leq s \leq N-k$. Then by using (26) we have $\mathbb{E}(V_{t+s}|\mathcal{F}_t) \leq K(1 + \rho^s |\mathcal{X}_{t,N}|_1)$. Therefore for any $R_1$ and $s > 0$, we have

$$\mathbb{E}(V_{t+s}|\mathcal{F}_t) \leq \left\{ K\rho^s + \frac{K}{R_1} I(|\mathcal{X}_{t,N}|_1 > R_1) \right\} |\mathcal{X}_{t,N}|_1 + K I(|\mathcal{X}_{t,N}|_1 \leq R_1).$$

Thus we have

$$\mathbb{E}(V_{t+s}|\mathcal{F}_t) \leq \left( K\rho^s + \frac{K}{R_1} \right) I(|\mathcal{X}_{t,N}|_1 > R_1) \} |\mathcal{X}_{t,N}|_1 + K(\rho^s + 1) I(|\mathcal{X}_{t,N}|_1 \leq R_1) |\mathcal{X}_{t,N}|_1.$$

By choosing an appropriate $R_1$ we can find an $s_1$ such that, for all $s \geq s_1$, we have $K\rho^s + K/R_1 < 1$ and thus condition (a) is satisfied. Condition (b) directly follows from Lemma A.3. Let $\mathbb{I}_{p+1}$ be a $(p+1) \times (p+1)$ matrix where $(\mathbb{I}_{p+1})_{ij} = 1$ for $1 \leq i, j \leq (p+1)$. Since $\|F_{l,N}\| = \mathcal{X}_{l,N}^\top \mathcal{X}_{l,N} / |\mathcal{X}_{l,N}|_1^2 = 1$, for $\lambda < 1$, we have $\lambda \|F_{l,N}\| < 1$, hence condition (c) is satisfied. Finally, by using the above for any $q \geq 1$ and $t \in \{N-p, \ldots, t\}$, we have

$$\sum_{l=t}^{r+s_1} \mathbb{E}(\|F_{l,N}\|^q | \mathcal{F}_j) \leq k_1 (p+1)^{2q}.$$

Therefore condition (d) is satisfied and Theorem 4.1 follows from Moulines *et al.* ([11], Theorem 16). $\square$

### A.3. The lower-order terms in the pertubation expansion

In this section we will prove the auxillary results required in Section 4.2, where we showed that the second-order terms in the pertubation expansion were a of lower order than the principle terms. The analysis of $J_{t_0,N}^{x,2}$ is based on partitioning it into two terms,

$$J_{t_0,N}^{x,2} = A_{t_0,N}^x + B_{t_0,N}^x,$$

similar to the partition in (52). $A_{t_0,N}^x$ is the weighted sum of $\{F_{k,N} - F_k(u)\}$, whereas $B_{t_0,N}^x$ is the weighted sum of the differences between the stationary approximation $F_k(u_0)$ and $\mathbb{E}(F_0(u_0))$, that is, of $\bar{F}_k(u_0) = F_k(u_0) - F(u_0)$. In this section we evaluate bounds for these two terms. We require the following lemma.

**Lemma A.4.** *Suppose Assumption 2.1 holds with some $r \geq 1$. Let $\mathcal{M}_{t,N}$ and $\mathcal{M}_t(u)$ be defined as in (29) and let $D_{t,N}$ be defined as in (45). Then we have*

$$\frac{\mathcal{X}_{t,N}\mathcal{X}_{t,N}^\top}{|\mathcal{X}_{t,N}|_1^2} = \frac{\mathcal{X}_t(u)\mathcal{X}_t(u)^\top}{|\mathcal{X}_t(u)|_1^2} + R_{t,N}(u) \tag{105}$$

*and*

$$\mathcal{M}_{t,N} = \mathcal{M}_t(u) + (1 + Z_t^2)R_{t,N}(u),$$

*where*

$$|R_{t,N}(u)| \leq \left( \left| \frac{t}{N} - u \right|^\beta + \left( \frac{p}{N} \right)^\beta + \frac{1}{N^\beta} \right) D_{t,N}.$$

*Further, for $q \leq q_0$ there exists a constant $K$ independent of $t, N$ and $u$ such that*

$$(\mathbb{E}|R_{t,N}(u)|^q)^{1/q} \leq K \left( \left| \frac{t}{N} - u \right|^\beta + \left( \frac{p+1}{N} \right)^\beta \right). \tag{106}$$

**Proof.** The proof uses (21) and the method given in Dahlhaus and Subba Rao ([4], Lemma A.4). □

We now give a bound for a general $A_{t,N}^x$.

**Lemma A.5.** *Suppose Assumption 2.1 holds with $r = q_0$ and let $\{G_{t,N}\}$ be a random process which satisfies $\sup_{t,N} \|G_{t,N}\|_{q_0}^E < \infty$. Let $F_{t,N}$, $F_t(u)$ and $F(u)$ be defined as in (30) and (9), respectively. Then if $|t_0/N - u_0| < 1/N$ and $q \leq q_0$ we have*

$$\left( \mathbb{E} \left| \sum_{k=0}^{t_0-p-1} (I - \lambda F(u_0))^k (F_{t_0-k-1,N} - F_{t_0-k-1}(u_0))G_{t_0-k,N} \right|^{q/2} \right)^{2/q}$$

$$\leq \frac{K}{(\delta N \lambda)^\beta} \sup_{t,N} (\mathbb{E}|G_{t,N}|^q)^{1/q}, \tag{107}$$

*where $K$ is a finite constant.*

**Proof.** By using (99) and (105) we obtain the result. □

In order to bound $\mathbb{E}|B_{t_0,N}^x|^q$ we need a Burkhölder type inequality (using Minkowski's inequality is not sufficient). This inequality is embedded in the following lemma which is a generalization of Proposition B.3 in Moulines *et al.* [11]. It can be proved by adapting the proof in Moulines *et al.* [11]; see also Dahlhaus and Subba Rao [3].

**Lemma A.6.** *Suppose $\{M_t\}$ and $\{F_t\}$ are random matrices and $F$ is a positive definite, deterministic matrix, with $\lambda_{\min}(F) > \delta$, for some $\delta > 0$. Let $\mathcal{F}_t = \sigma(F_t, M_t, F_{t-1}, M_{t-1}, \ldots)$. Assume, for some $q \geq 2$, the following:*

  (i) *$\{F_t\}$ are identically distributed with mean zero.*
 (ii) *$\mathbb{E}(M_t | \mathcal{F}_{t-1}) = 0$.*
(iii) *$(\mathbb{E}|\mathbb{E}(F_t|\mathcal{F}_k)|^{2q})^{1/2q} \leq K\rho^{t-k}$*
 (iv) *$(\mathbb{E}|\mathbb{E}(M_t F_s | \mathcal{F}_k) - \mathbb{E}(M_t F_s)|^q)^{1/q} \leq K\rho^{s-k}$ if $k \leq s \leq t$.*
  (v) *$\sup_t(\mathbb{E}|M_t|^{2q})^{1/(2q)} < \infty$ and $(\mathbb{E}|F_t|^{2q})^{1/(2q)} < \infty$.*

*Then we have*

$$\left( \mathbb{E} \left| \sum_{k=0}^{t-p-1} \sum_{i=0}^{t-p-k-2} (I - \lambda F)^{k+i} F_{t-k-1} M_{t-k-i-1} \right|^q \right)^{1/q} \leq \frac{K}{\delta\lambda}, \tag{108}$$

*where $K$ is a finite constant.*

We now apply the lemma above to the particular example of the ANRE algorithm.

**Lemma A.7.** *Suppose Assumption 2.1 holds for $r > 4$. Let $\{F_t(u)\}$, $\{\mathcal{M}_{t,N}\}$ and $\{\mathcal{M}_t(u)\}$ be defined as in (29) and $\bar{F}_t(u) = F_t(u) - \mathbb{E}(F_t(u))$. Then*

$$\left( \mathbb{E} \left| \sum_{k=0}^{t_0-p-1} \sum_{i=0}^{t_0-p-k-2} \{I - \lambda F(u)\}^{k+i} \bar{F}_{t_0-k-1}(u) \mathcal{M}_{t_0-k-i-1,N} \right|^{r/2} \right)^{2/r} \leq \frac{K}{\lambda} \tag{109}$$

*and*

$$\left( \mathbb{E} \left| \sum_{k=0}^{t_0-p-1} \sum_{i=0}^{t_0-p-k-2} \{I - \lambda F(u)\}^{k+i} \bar{F}_{t_0-k-1}(u) \mathcal{M}_{t_0-k-i-1}(u) \right|^{r/2} \right)^{2/r} \leq \frac{K}{\lambda}. \tag{110}$$

**Proof.** We prove (109); the proof of (110) is the same. We verify the conditions in Lemma A.6, then (109) immediately follows. By using (97) we have that $\lambda_{\min}\{F(u)\} \geq \delta$, for some $\delta > 0$. Let $M_t := \mathcal{M}_{t,N}$, $F_t := \bar{F}_t(u)$, $F := F(u)$ and $\mathcal{F}_t = \sigma(Z_t, Z_{t-1}, \ldots)$. It is clear from the definition that the $\{\bar{F}_t(u)\}_t$ have zero mean and are identically distributed; also $\mathbb{E}(\mathcal{M}_t(u)|\mathcal{F}_{t-1}) = 0_{p+1 \times p+1}$. By using (24) we have

$$\left(\mathbb{E}|\mathbb{E}(\bar{F}_t(u)|\mathcal{F}_k)|^r\right)^{1/r} = \left(\mathbb{E}|\mathbb{E}(F_t(u)|\mathcal{F}_k) - \mathbb{E}(F_t(u))|^r\right)^{1/r} \leq K\rho^{t-k},$$

thus condition (iii) is satisfied. Since $\mathcal{M}_{t,N} = (Z_t^2 - 1)\sigma_{t,N}^2 \mathcal{X}_{t-1,N}/|\mathcal{X}_{t-1,N}|_1^2$, and $F_t \leq \mathbb{I}_{p+1}$ and $\sigma_{t,N}^2 \mathcal{X}_{t-1,N}/|\mathcal{X}_{t-1,N}|_1^2 \leq \mathbb{I}_{p+1}$, by using Corollary A.1 and $\sup_{t,N}(\mathbb{E}|\mathcal{X}_{t,N}|^{r/2})^{2/r} < \infty$, we can show that condition (iv) is satisfied. Moreover, $F_{t,N}$ is a bounded random matrix, hence all its moments exist. Finally, since $\sup_{t,N}|\mathcal{M}_{t,N}|_r^r \leq K(Z_t^2 + 1)^r$ and $\mathbb{E}(Z_0^{2r}) < \infty$, we have, for all $k \leq s \leq t \leq N$, that $\sup_{t,N}(\mathbb{E}|\mathcal{M}_{t,N}|^r) < \infty$, leading to condition (v). Thus all the conditions of Lemma A.6 are satisfied and we obtain (109). $\square$

# Acknowledgements

This work has been supported by the Deutsche Forschungsgemeinschaft (DA 187/12-2). The authors would like to thank two anonymous referees for make several constructive comments which improved the presentation of the paper.